\documentclass{amsart}
\usepackage{amssymb,latexsym}

\newtheorem{theorem}{Theorem}[section]
\newtheorem{proposition}[theorem]{Proposition}
\newtheorem{corollary}[theorem]{Corollary}
\newtheorem{definition}[theorem]{Definition}

\newtheorem{remark}[theorem]{Remark}
\newtheorem{lemma}[theorem]{Lemma}
\newtheorem{example}[theorem]{Example}

\newcommand{\torus}{{\tilde{\mathcal C}}/{\mathbb Z}^d}

\def\proof{\smallskip\noindent {\it Proof: \ }}
\def\proofof#1{\smallskip\noindent {\it Proof of #1: \ }}
\def\endproof{\hfill$\square$\medskip}

\def\A{\mathcal{A}}
\def\C{\mathcal{C}}
\def\T{\mathcal{T}}
\def\M{\mathcal{M}}

\def\x{\mathbf{x}}
\def\y{\mathbf{y}}
\def\k{\mathbf{k}}

\def\P{\mathbb{P}}
\def\R{\mathbb{R}}
\def\T{\mathbb{T}}
\def\Z{\mathbb{Z}}
\def\N{\mathbb{N}}
\def\ZZ{\mathcal{Z}} 

\def\rank{\mathrm{rk}}

\def\O{\Delta}
\def\<{\left<}
\def\>{\right>}

\def\lcm{\mbox{\upshape {lcm}}\,}
\def\field{{\bf k}}
\def\ord{\mbox{\upshape {\small {ord}}}\,}

\def\depth{\mbox{\upshape {depth}}\,}

\def\matroid{\underline{\mathcal{M}}}
\def\L{\mathcal{L}}

\title{Syzygies of Oriented Matroids}

\author{Isabella Novik 
\and    Alexander Postnikov
\and    Bernd Sturmfels}

\address{Department of Mathematics,
         University of California,
         Berkeley, CA 94720}

\email{bellan@math.berkeley.edu}
\email{apost@math.berkeley.edu}
\email{bernd@math.berkeley.edu}

\begin{document}

\begin{abstract}
We construct minimal cellular resolutions of squarefree monomial ideals 
arising from hyperplane arrangements, matroids and oriented matroids.
These are Stanley-Reisner ideals of complexes of independent sets,
and of triangulations of Lawrence matroid polytopes.
Our resolution provides a cellular realization of  Stanley's formula 
for their Betti numbers. For unimodular matroids our resolutions are
related to hyperplane arrangements on tori, and we recover the resolutions
constructed by Bayer, Popescu and Sturmfels \cite{BPS}.
We resolve the combinatorial problems posed in  \cite{BPS}
by computing M\"obius invariants 
of graphic and cographic arrangements in terms of
Hermite polynomials.
\end{abstract}

\maketitle

\section{Cellular Resolutions from Hyperplane Arrangements}
\label{sec:monomial}

A basic problem of combinatorial commutative algebra is to find the
syzygies of a monomial ideal $\,M \,= \langle m_1,\ldots,m_r \rangle \,$
in the polynomial ring $\k[\x]=\k[x_1,\dots,x_n]$ over a field $\k$. 
One approach involves constructing {\it cellular resolutions},
where the $i$-th syzygies of $M$ are indexed by the $i$-dimensional
faces of a CW-complex on $r$ vertices. After reviewing
the general construction of cellular resolutions from \cite{bs},
we shall define the monomial ideals and resolutions that are studied
in this paper.

Let~$\Delta$ be a {\it CW-complex} \cite[\S 38]{M}
with $r$ vertices   $v_1,\dots,v_r$, which are labeled by the monomials
$m_1, \ldots, m_r$.  
We write $c\geq c'$ whenever a cell~$c'$ belongs to the closure of
another cell $c$ of $\Delta$. This defines the face poset of $\Delta$.
We label each cell $c$ of $\Delta$ with 
the monomial $m_c=\lcm\{m_i\mid v_i\leq c\}$,  the
least common multiple of the monomials labeling the vertices of~$c$. 
Also set $m_\emptyset =1$ for the empty cell of $\Delta$.  
Clearly, $m_{c'}$ divides $m_c$ whenever $c'\leq c$.
The principal ideal $\<m_c\>$ is identified with the
free $\N^n$-graded $\k[\x]$-module of rank $1$ with generator in degree
$\deg m_c$.  For a pair of  cells
$c\geq c'$, let $p_c^{c'}:\<m_c\>\to \<m_{c'}\>$ 
be the inclusion map of ideals.  
It is a degree-preserving homomorphism of $\N^n$-graded modules.

Fix an orientation of each cell in $\Delta$,
and define the {\it cellular complex\/}  $\, C_\bullet(\Delta,M) \,\, $
$$ \cdots \,{\overset{\partial_3}\longrightarrow}\, C_2
\,{\overset{\partial_2}\longrightarrow}\, C_1
\,{\overset{\partial_1}\longrightarrow}\, C_0
\,{\overset{\partial_0}\longrightarrow}\, C_{-1}
=\k[\x]
$$
as follows.
The $\N^n$-graded $\k[\x]$-module of $i$-chains is 
$$
C_i \quad = \,\,\, \bigoplus_{c\,:\, \dim c = i} \<m_c\>\,,
$$
where the direct sum is over all $i$-dimensional cells~$c$ of~$\Delta$.
The differential $\partial_i:C_{i}\to C_{i-1}$ is defined 
on the component $\<m_c\>$ as the alternating sum of
the maps ~$p_c^{c'}$: 
$$
\partial_i \quad = \sum_{c'\leq c,\,\dim c' =i-1} 
[c:c'] \ p_c^{c'}\,,
$$
where $[c:c']\in\Z$ is the {\it incidence coefficient\/} of oriented
cells $c$ and $c'$ in the usual topological sense.
For a regular CW-complex, the incidence coefficient $[c:c']$
is $+1$ or $-1$ depending on the orientation of cell $c'$ in the boundary
of $c$.  The differential $\partial_i$ preserves the $\N^n$-grading
of $\k[\x]$-modules.
Note that if $m_1=\cdots=m_r=1$ then $C_\bullet(\Delta,M)$ is the usual 
chain complex of~$\Delta$ over $\k[\x]$.
For any monomial $m \in \k[\x]$, we define
$\Delta_{\leq m}$ to be the subcomplex of $\Delta $ consisting
of all cells $c$ whose label $m_c$ divides $m$.
We call any such $\Delta_{\leq m}$ an
{\it $M$-essential} subcomplex of $\Delta$.

\begin{proposition}  {\rm \cite[Proposition~1.2]{bs}} \
\label{prop:resolution}
The cellular complex $C_\bullet(\Delta,M)$ is exact
if and only if every $M$-essential subcomplex
$\Delta_{\leq m}$ of $\Delta$ is  acyclic over~$\k$.
Moreover, if $m_c\ne m_{c'}$ for any $c>c'$, then 
$C_\bullet(\Delta,M)$ gives a minimal free resolution of $M$.
\end{proposition}  

Proposition~\ref{prop:resolution} is derived from
the observation that, for a monomial $m$,
the $(\deg m)$-graded component 
of~$C_\bullet(\Delta,M)$ equals the chain complex 
of $\Delta_{\leq m}$ over $\k$. If both of the hypotheses in
Proposition~\ref{prop:resolution} are met, then we say that
$\Delta$ is an {\it $M$-complex}, and we call
$C_\bullet(\Delta,M)$ a {\it minimal cellular resolution\/} of $M$.
Thus each $M$-complex $\Delta$ produces 
a minimal free resolution of the ideal~$M$.
In particular, for an $M$-complex $\Delta$, the number 
$f_i(\Delta)$ of $i$-dimensional cells of $\Delta$ is exactly the {\it $i$-th 
Betti number\/} of~$M$, i.e., the rank of the $i$-th free module in 
a minimal free resolution.  Thus, for fixed $M$, all $M$-complexes
have the same $f$-vector.

Examples of $M$-complexes appearing in the literature include
planar maps \cite{MS},  Scarf complexes \cite{Peeva}
 and hull complexes \cite{bs}.
A general construction of $M$-complexes using  discrete Morse theory
was proposed by Batzies and Welker \cite{BW}.
We next introduce a family of $M$-complexes 
which generalizes those in \cite[Theorem~4.4]{BPS}.

Let~$\A = \{H_1,H_2,\ldots,H_n\}$ be an arrangement of 
$n$ affine hyperplanes in $\R^d$, 
\begin{equation}
\label{eq:affine-arrangement}
  H_i \quad = \quad \{v\in  \R^d \mid h_i (v) =c_i\},\quad i=1,\dots,n,
\end{equation}
where $c_1,\dots,c_n\in\R $ and
$h_1,\dots,h_n$ are nonzero linear forms that span~$(\R^d)^*$.
%, and the intersection $\,H_1 \cap \cdots \cap H_n\,$ of all
%hyperplanes is empty.

We fix two sets of variables $x_1,\dots,x_n$ and $y_1,\dots,y_n$, 
and we associate with the arrangement~$\A$ two functions
$m_x$ and $m_{xy}$ from $\R^d$ to sets of monomials:
$$
%\begin{array}{l}
 m_x:v\longmapsto \prod_{i\,:\,h_i (v)\ne c_i} \!\! x_i
\quad\textrm{and}\quad
%\\[.3in]
 m_{xy}:v\longmapsto \left(\prod_{i\,:\, h_i (v) > c_i} \!\! x_i \right) \cdot
 \left(\prod_{j\,:\, v_j(v) < c_j} \!\! y_j \right)\,\\[.01in].
%\end{array}
$$
Note that  $m_x(v)$ is obtained  from  $m_{xy}(v)$
by specializing $y_i$ to $x_i$ for all $i$.

\begin{definition} \rm
The {\it matroid ideal} of $\A$ is the ideal
$M_\A$ of $\k[\x]=\k[x_1,\dots,x_n]$ generated by the
monomials $\{m_x(v) : v\in \R^d \}$.
The {\it oriented matroid ideal} of $\A$ is the ideal
$O_\A$ of $\k[\x,\y]=\k[x_1,\dots,x_n,y_1,\dots,y_n]$ 
generated by $\{m_{xy}(v) : v\in \R^d \}$.
\end{definition}

The hyperplanes $H_1, \ldots, H_n$ partition $\R^d$ into relatively open 
convex polyhedra, called the {\it cells} of $\A$. Two points
$v,v'\in \R^d$ lie in the same cell $c$ if and only if $m_{xy}(v)=m_{xy}(v')$.
We write $m_{xy}(c)$ for that monomial, and similarly $m_x(c)$
for its image under $y_i \mapsto x_i$.
Note that $m_x(c')$ divides $m_{x}(c)$, 
and $m_{xy}(c')$ divides $m_{xy}(c)$, provided $c'\leq c$.
The cells of dimension $0$ and $d$ are called
{\it vertices\/} and {\it regions,} respectively.  A cell is
{\it bounded\/} if it is bounded as a subset of $\R^d$.
The set of all bounded cells forms a regular CW-complex $B_\A$
called the {\it bounded
complex\/} of $\A$.

\setlength{\unitlength}{1.5pt}
\begin{figure}[htbp]
\begin{center}
    \begin{picture}(100,195)(0,-55)

     \thinlines
      \put(0,0){\circle*{4}}
      \put(3,-8){$y_4$}

      \put(100,0){\circle*{4}}
      \put(90,-8){$y_2 x_3$}

      \put(0,100){\circle*{4}}
      \put(-17,95){$x_1 x_2$}
      
      \put(50,50){\circle*{4}}
      \put(55,48){$x_1 x_3$}

      \put(-40,0){\line(1,0){70}}
      \put(70,0){\line(1,0){70}}
      \put(-30,0){\vector(0,1){10}}
      \put(-55,0){$H_1$}

      \put(-30,-30){\line(1,1){50}}
      \put(30,30){\line(1,1){50}}
      \put(-20,-20){\vector(-1,1){7}}
      \put(-43,-40){$H_2$}

      \put(0,-40){\line(0,1){40}}
      \put(0,100){\line(0,1){40}}
      \put(0,-30){\vector(1,0){10}}
      \put(-4,-53){$H_3$}

      \put(-30,130){\line(1,-1){30}}
      \put(100,0){\line(1,-1){30}}
      \put(120,-20){\vector(1,1){7}}
      \put(133,-40){$H_4$}

      \put(41,-1){$y_2x_3y_4$}
      \put(-12,44){$x_1x_2y_4$}
      \put(18,74){$x_1x_2x_3$}
      \put(65,24){$x_1y_2x_3$}
      \put(16,24){$x_1x_3y_4$}

      \put(7,56){$x_1x_2x_3y_4$}
      \put(38,14){$x_1y_2x_3y_4$}

      \put(50,120){$M_\A=\left<x_4,x_1x_2,x_2x_3,x_1x_3\right>$}
      \put(50,100){$O_\A=\left<y_4,x_1x_2,y_2x_3,x_1x_3\right>$}

      \thicklines
      \put(0,0){\line(1,0){37}}
      \put(63,0){\line(1,0){37}}
      \put(0,0){\line(0,1){40}}
      \put(0,50){\line(0,1){50}}
      \put(0,0){\line(1,1){20}}
      \put(30,30){\line(1,1){20}}
      \put(0,100){\line(1,-1){20}}
      \put(30,70){\line(1,-1){40}}
      \put(80,20){\line(1,-1){20}}

    \end{picture}
    \caption{The bounded complex $B_\A$ with monomial labels.}
    \label{fig:arrangement}
\end{center}
\end{figure}
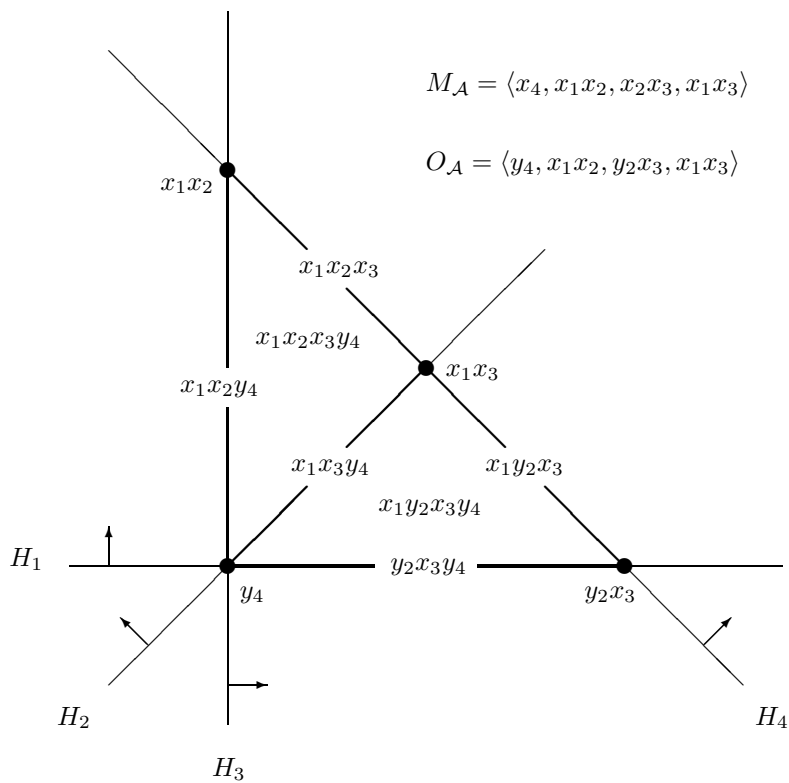

Figure~\ref{fig:arrangement} shows an example of a hyperplane 
arrangement~$\A$ with $d=2$ and $n=4$, 
together with monomials that label its bounded cells. 
The bounded complex $B_\A$ of this arrangement consists of 4 vertices, 
5 edges, and 2 regions.

\begin{theorem} \
\label{th:arrang_resolution}
{\rm\bf (a)} \ The ideal $M_\A$ is minimally generated by the monomials
$m_x(v)$, where $v$ ranges over the vertices of~$\A$.  
The bounded complex $B_\A$ is an $M_\A$-complex.
Thus its cellular complex $C_\bullet(B_\A,M_\A)$
gives a minimal free resolution for $M_\A$.

\medskip
\noindent
{\rm\bf (b)} \ The ideal $O_\A$ is minimally generated by the  monomials
$m_{xy}(v)$, where $v$ ranges over the vertices of~$\A$.  
The bounded complex $B_\A$ is
an $O_\A$-complex.  Thus its cellular complex $C_\bullet(B_\A,O_\A)$
gives a minimal free resolution for $O_\A$.
\end{theorem} 

To prove Theorem~\ref{th:arrang_resolution}, we must check
that for both ideals, the two hypotheses of 
Proposition~\ref{prop:resolution} are satisfied.
The second hypothesis is immediate:
for a pair of cells $c>c'$, there is a hyperplane~$H_i \in \A$
that contains $c'$ but does not contain $c$, in which case
$m_x(c)$ is divisible by $x_i$ and  $m_x(c')$ is not divisible by $x_i$.
Analogously, for the oriented matroid ideal~$O_\A$.
The essence of Theorem~\ref{th:arrang_resolution} is the acyclicity 
condition, which states that all $M_\A$-essential 
and $O_\A$-essential subcomplexes of $B_\A$ are acyclic.
For the whole bounded complex this is known:

\begin{proposition}  \label{prop:contractible}
\label{prop:B-contractible}
{\rm (Bj\"{o}rner and Ziegler,
see \cite[Theorem~4.5.7]{BLSWZ})} \hfill \break
The complex $B_\A$ of bounded cells of a hyperplane 
arrangement $\A$ is contractible.
\end{proposition}

The acyclicity of all $M_\A$-essential subcomplexes of $B_\A$
is an easy consequence of Proposition~\ref{prop:contractible}:
each $M_\A$-essential subcomplex is a bounded complex of a
hyperplane arrangement induced by $\A$ in one of the flats of $\A$.
The acyclicity of all $O_\A$-essential subcomplexes will follow from
a  generalization of Proposition~\ref{prop:contractible}
stated in Proposition~\ref{contractible} below.
We give more details in Section~\ref{oriented_section}, where 
Theorem~\ref{th:arrang_resolution} is restated and proved
in the more general setting of oriented matroids.

\vskip .1cm

The main result in this paper is the
construction of the minimal free resolution of 
an arbitrary matroid ideal (Theorems \ref{res2} and \ref{S-exact})
and an arbitrary oriented matroid ideal (Theorem \ref{res1}).
A numerical consequence of this result is a refinement of
Stanley's formula, given in \cite[Theorem~9]{St77}, 
 for their Betti numbers (Corollary~\ref{O-Hilbert},
Corollary~\ref{happiness}; see also the last paragraph of Section~3).
The simplicial complexes corresponding to matroid ideals and
oriented matroid ideals 
are the complexes of independent sets in matroids
(Remark \ref{easyremark}) and the triangulations of 
Lawrence matroid polytopes (Theorem \ref{lawlaw}), respectively.
In the unimodular case, oriented matroid ideals 
arise as initial ideals of toric varieties
in $\P^1 \times \P^1 \times \cdots \times \P^1$, by
work of Bayer, Popescu and Sturmfels 
\cite[\S 4]{BPS}, and their Betti numbers can be interpreted 
as face numbers of hyperplane arrangements on a torus
(Theorem \ref{f-vector}). Every ideal considered in this
paper is Cohen-Macaulay, its Cohen-Macaulay type
(= highest Betti number) is the M\"obius invariant of 
the underlying matroid, and all other Betti numbers
are sums of M\"obius invariants of matroid minors
(\S 4, (\ref{stanley})). In Section 5 we resolve
the enumerative problems concerning graphic and cographic 
matroids which were left open in \cite[\S 5]{BPS}.
Propositions  \ref{prop:only-source} and  \ref{Coinvariant}
give combinatorial expressions for the M\"obius invariant of any graph.
More precise and efficient formulas, in terms of Hermite polynomials,
are established for the M\"obius coinvariants of complete graphs
(Theorem \ref{th:mu(K_m)}) and of complete bipartite graphs
(Theorem \ref{th:mu-bipartite}).

\bigskip

\section{Oriented Matroid Ideals} 
\label{oriented_section} 

The axiomatic theory of oriented matroids provides a foundation 
for geometric combinatorics, much in the same way as the
axiomatic theory of schemes provides a foundation for algebraic geometry.
Oriented matroid techniques are ubiquitous in the study of hyperplane 
arrangements, point configurations, and convex polytopes. 
In this section we establish a link between  oriented matroids and
commutative algebra. In the resulting combinatorial context,
the algebraists' classical question, ``What makes a complex exact?'' 
\cite{BE},  receives a surprising answer: it is
the Topological Representation 
Theorem of Folkman and Lawrence \cite[Chapter 5]{BLSWZ}.

We  start by briefly reviewing one of the axiom systems for 
oriented matroids~\cite{BLSWZ}. 
Fix a finite set $E$. A {\it sign vector} $X$ is an element of
$\{+,-,0\}^E$. The {\it positive part}
 of $X$ is denoted $X^+=\{i\in E : X_i=+\}$,
and similarly $X^-$ and $X^0$.  The support of $X$ is
$\underline{X}=\{i\in E : X_i\neq 0\}$.
The {\it opposite} $-X$ of a vector $X$ is given by $(-X)_i=-X_i$.
The {\it composition} $X \circ Y$ of two vectors $X$ and $Y$ is the sign
vector defined by
\[
%(-X)_i =\left\{ \begin{array}{lll}
%              -    &\mbox{if $X_i =+$} \\
%              +    &\mbox{if $X_i=-$} \\
%              0    &\mbox{if $X_i=0$}
%            \end{array}
%    \right.  \quad
%\mbox{ and }
(X\circ Y)_i =\left\{ \begin{array}{ll}
              X_i    &\textrm{if $X_i \neq 0$,} \\
              Y_i    &\textrm{if $X_i=0$.} 
            \end{array}
    \right.  
%\mbox{ respectively}.
\]
The {\it separation set\/} of sign vectors $X$ and $Y$ is
$S(X,Y)=\{i\in E\mid X_i=-Y_i\ne 0\}$.  

A set $\L \subseteq \{+,-,0\}^E$ is the set of {\em covectors}
of an {\it oriented matroid on $E$}  if and only if it 
satisfies the following four axioms \cite[\S~4.1.1]{BLSWZ}:
\begin{enumerate}
\item the zero sign vector $0$ is in $\L$;
\item (symmetry) \ if $X\in\L$ then $-X\in\L$;
\item (composition) \ if $X,Y\in\L$ then $X\circ Y\in\L$;
\item (elimination) \ if $ X,Y\in\L$ and $i\in S(X,Y)$
then there exists $Z\in\L$ such that $Z_i=0$ and
$Z_j=(X\circ Y)_j=(Y\circ X)_j$ for all $j\not\in S(X,Y)$.
\end{enumerate}
Somewhat informally, we  say that such a pair $(E, \L)$,
{\it is} an oriented matroid. 
An {\it affine oriented matroid\/} \cite[\S 10.1]{BLSWZ}, 
denoted $\M=(E,\L,g)$, is an oriented matroid 
with a distinguished element $g\in E$
such that $g$ is not a {\it loop,} i.e., 
$X_g \ne 0$ for at least one covector $X\in\L$.
The {\it positive part\/} of $\L$  is  
$\L^+=\{X \in \L \; : \: X_g=+\}$.

The set $\{+, -, 0\}^E$ is partially ordered by the product of partial
orders
$$
              0 < + \quad\textrm{and}\quad
              0 < -     \qquad 
              \textrm{($+$ and $-$ are incomparable)}.
$$
This induces a partial order on the set of covectors $\L$.
A covector $X$ is called {\it bounded\/} if every nonzero covector
$Y\leq X$ is in the positive part $\L^+$.

The Topological Representation Theorem for Oriented Matroids, see
\cite[Theorem~5.2.1]{BLSWZ}, states that  $\widehat\L=\L\cup\{\hat 1\}$
is the face lattice of an arrangement of pseudospheres,
and  $\widehat\L^+=\L^+\cup\{\hat 0,\hat 1\}$ 
is the face lattice of an arrangement of pseudohyperplanes
\cite[Exercise 5.8]{BLSWZ}. These are 
regular CW-complexes homeomorphic to a sphere and a ball, respectively.
(This is why $\widehat\L$ is called the {\it face lattice},
and $\widehat\L^+$ is called the  {\it affine face lattice} of $\M$).
The {\it bounded complex} $B_\M$ of $\M$ is their subcomplex
 formed by the cells associated with the bounded covectors.
The bounded complex is uniquely determined by its face lattice---the
poset of bounded covectors.   Slightly abusing notation, we denote
this poset by the same symbol $B_\M$.

%The lattice  $\widehat\L$
%and its sublattice $\widehat\L^+$ is graded. 

We write $\rank(\,\cdot \,)$ for the
rank function of the lattice  $\widehat\L$. The atoms of $\widehat\L$, i.e.,
the elements of rank $1$, are called {\it cocircuits\/}
of~$\M$.  The vertices of the bounded complex $B_\M$
are exactly the cocircuits of $\M$ that belong to the positive part $\L^+$.

\begin{example}    
\label{central}
{\rm  (Affine oriented matroids  from hyperplane arrangements) \
\hfill \break
Let  $\C=\{H_1, \ldots, H_n, H_g\}$ be a central hyperplane arrangement
in $\R^{d+1} = \R^d \times \R$, written as
$H_i =\{(v,w) \in \R^d\times \R : h_i (v) =c_i  w\}\,$ 
and $H_g = \{(v,w) : w = 0 \}$.
The restriction of $\C$ to the hyperplane $\{(v,w)  : w = 1\}$
is precisely  the affine arrangement $\A$ in Section 1.
Fix $E = \{1,\ldots,n,g \}$. The image of the map
$$  \R^{d+1} \rightarrow \{+,-,0\}^E \, , \,\,
(v,w) \mapsto \bigl(\,
{\rm sign}(h_1(v)  -  c_1 w),\, \ldots,\,
{\rm sign}(h_n(v) -  c_n w), \,{\rm sign}(w) \,\bigr) $$
is the set $\L$ of covectors of an oriented matroid on $E$.
The affine face lattice $\widehat{\L}^+$ of
$\M=(E, \L, g)$ equals the face lattice of
the affine hyperplane arrangement $\A$. The
 bounded complex $B_\M$
coincides with the bounded complex $B_\A$ in
Proposition \ref{prop:B-contractible}.
}
\end{example}

Let  $\M=(E, \L, g)$  be an
affine oriented matroid on 
$E=\{1, \ldots, n, g\}$. 
With every sign vector $Z\in \{0,+, -\}^E$ we associate a monomial 
$$
m_{xy}(Z) = \left(\prod_{\,i\,:\,Z_i=+} x_i \right)\cdot
\left(\prod_{\,j\,:\,Z_i=-} y_i\right)\,,
\qquad \quad \hbox{where  $x_g=y_g=1$.}
$$
The {\it oriented matroid ideal\/} $O$ is
the ideal in the polynomial ring 
$\k[\x,\y]=\field[x_1,\dots,$ $ x_n, y_1, \dots, y_n]$ 
generated by all monomials corresponding to
covectors $Z\in\L^+$.
The {\it matroid ideal\/} $M$ associated with 
$\M=(E, \L, g)$ is the ideal of $\k[\x]$ obtained from
$O$ by specializing $y_i$ to $x_i$ for all~$i$.
These ideals are treated in Section~\ref{matroid_section}.
The main result of this section concerns the syzygies of
the oriented matroid ideal $O$.

\begin{theorem}  \label{res1}
The oriented matroid ideal $O$ 
is minimally generated by the monomials corresponding to
the vertices of $B_\M$. The bounded complex
$B_\M$ is an $O$-complex. Thus its cellular complex
$C_\bullet(B_\M, O)$ gives a minimal 
$\N^{2n}$-graded free $\k[\x,\y]$-resolution of~$O$.

\end{theorem}

Recall that for a monomial $m$ in $\k[\x,\y]$, the corresponding 
$\N^{2n}$-graded Betti number of $O$, $\beta_m(O)$, is the 
multiplicity of the summand $\langle m \rangle$ in a minimal 
$\N^{2n}$-graded $\k[\x,\y]$-resolution of~$O$.
Theorem~\ref{res1} implies the following numerical result.

\begin{corollary}  \label{O-Hilbert}
The $\N^{2n}$-graded Betti numbers of $O$ are all $0$ or $1$.  
They are given by 
%the absolute values of 
the coefficients in the numerator of the $\N^{2n}$-graded 
Hilbert series of~$O$:
\begin{equation}
\label{hilbertseries}
 \biggl(\sum_{Z\in B_\M} (-1)^{\rank(Z)}
m_{xy}(Z) \biggr) / \prod_{i=1}^n (1-x_i)(1-y_i).
\end{equation}
\end{corollary}

{\smallskip\noindent {\it Proof of Theorem \ref{res1}: \ }}
Distinct cells  $Z$ and $Z'$ of the bounded complex $B_\M$ 
have  distinct labels $m_{xy}(Z) \neq m_{xy}(Z')$. 
% For $Z\geq Z'$, the monomial $m_{xy}(Z)$ is divisible by  $m_{xy}(Z')$.
This implies minimality of the complex $C_\bullet(B_\M, O)$. 
%We conclude that all $\N^{2n}$-graded Betti numbers are $0$ or $1$, 
%and we obtain the asserted Hilbert series for $\k[\x,\y]/O$. 
In order to prove exactness of $C_\bullet(B_\M, O)$, we must verify the
first hypothesis in Proposition~\ref{prop:resolution}. 
To this end, we shall digress and first present a
generalization of Proposition~\ref{prop:contractible}. 

The {\it regions\/} of an oriented matroid
$(E, \L)$ are the maximal covectors, i.e, 
the maximal elements of 
the poset $\L$. For a covector $X\in\L$ and a subset
$E'$ of~$E$, denote by $ X|_{E'}\in \{+,-,0\}^{E'}$ the restriction 
of $X$ to $E'$: $(X|_{E'})_i = X_i$, for every $i\in E'$. 
The restriction of $(E, \L)$ to
a subset $E'$ of $E$ is the oriented matroid on $E'$ with the set
of covectors $\,\L|_{E'}= \{X|_{E'} : X\in\L\}$.

The following result,
which was cited without proof in \cite[Theorem~4.4]{BPS},
is implicit in the derivation  of \cite[Theorem~4.5.7]{BLSWZ}.
We are grateful to G\"unter Ziegler for making this explicit
by showing us the following proof. Ziegler's proof does not 
rely on the Topological Representation Theorem for Oriented Matroids.
If one uses that theorem, then 
Proposition \ref{contractible} can also be proved 
by a topological argument.

\begin{proposition} {\rm (G.~Ziegler)}  \label{contractible}
Let $\M=(E, \L, g)$ be an affine oriented matroid and
$B_\M$ its bounded complex. 
For any subset $E'$ of $E$ and any region $R'$ of $(E', \L|_{E'})$,
the CW-complex with the face poset
 $\,B'=\{X\in B_\M :
 X|_{E'} \leq R'\}\,$ is contractible. 
\end{proposition}  

\proof
Let $\T$ denote the set of regions of $\L$.
A subset $A \subseteq \T$ is said to be {\it $T$-convex\/}
if it is an intersection of ``half-spaces'', i.e., sets of the form
$\T^+_e = \{ T\in \T : T_e=+\}$ and $\T^-_e = \{ T\in \T : T_e=-\}$.
Each region $R\in \T$ defines a partial order on  $\T$:
$$   T_1 \le T_2  \quad : \iff \quad
  \{e\in E : R_e = -(T_1)_e \} \,\subseteq \,\{e\in E : R_e = -(T_2)_e \}.$$
Denote this poset by $\T(\L,R)$. We also abbreviate
$\T^+ := \T^+_g = \T\cap \L^+$.

We may assume that $B'$ is non-empty. Then 
$\mathcal{R}:=\{X\in\T^+: X|_{E'} =R'\}$ is a
non-empty, $T$-convex set. 
Lemma 4.5.5 in \cite{BLSWZ} states that
$\mathcal{R}$ is an order ideal of
$\T(\L,R)$, and, moreover, it is an order ideal of 
$\T^+ \subseteq\T(\L,R)$.
By \cite[Proposition~4.5.6]{BLSWZ}, 
there exists a recursive coatom ordering of $\widehat{\L}^+$
in which the elements of $\mathcal{R}$ come first.
The restriction of this ordering to $\mathcal{R}$ is a
recursive coatom ordering of the poset 
$ \widehat{\L}^+_{\mathcal{R}} =
 \{ X \in \L^+ :  X\leq T \mbox{ for some } T\in {\mathcal{R}} \} 
   \cup \{\hat{1}\}$.
This implies (using \cite[Lemma 4.7.18]{BLSWZ}) that the order complex 
$\Delta_{\ord}(\L^+_{\mathcal{R}})$ of  $\L^+_{\mathcal{R}}$
is a shellable $(r-1)$-ball. It is a subcomplex
of $\Delta_{\ord}(\L^+)$, which is also an $(r-1)$-ball,
by \cite[Theorem~4.5.7]{BLSWZ}. Let 
$U=\L^+_{\mathcal{R}} \backslash B_\M$ be the set
of ``unbounded covectors''. Then the subcomplex $\Delta_U$
of $\Delta_{\ord}( \L^+_{\mathcal{R}})$
induced on the vertex set of  $U$
 lies in the boundary of $\Delta_{\ord}(\L^+)$, and hence also
in the boundary of $\Delta_{\ord}(\L^+_{\mathcal{R}})$. 
Thus 
$||\Delta_{\ord}(\L^+_{\mathcal{R}})|| \setminus ||\Delta_U||$
is a contractible space. By \cite[Lemma 4.7.27]{BLSWZ}, the space
$||\Delta_{\ord}(B')||$ is a strong deformation retract of 
$||\Delta_{\ord}(\L^+_{\mathcal{R}})|| \backslash ||\Delta_U||$,
and is hence contractible as well.
\endproof

\vskip .2cm

We now finish the proof of Theorem~\ref{res1}. Consider any $O$-essential 
subcomplex $(B_\M)_{\leq \x^{\bf a}\y^{\bf b}}$
of $B_\M$, with
${\bf a}, {\bf b} \in \N^n$.
This complex consists of all cells $Z$
 whose label $m_{xy}(Z)$ divides $\x^\mathbf{a}\,\y^\mathbf{b}$.
Set 
\begin{eqnarray*}
E''&=&\{1 \leq i \leq n : a_i = 0 \mbox{ and } b_i = 0\},\\
E'&=& \{1 \leq i \leq n : \mbox{ exactly one of } a_i  \mbox{ and } 
   b_i \mbox{ is positive} \} \subseteq E\setminus E''.
\end{eqnarray*}
We first replace our affine oriented matroid
 $(E,\L,g)$ by the affine oriented matroid
$\,(E \backslash E'', \L / E'', g)\,$  gotten by contraction at $E''$.
Next we define $R'\in \{+,-,0\}^{E'}$ by
\[
R'_i = \left\{ \begin{array}{ll}
              +    &\mbox{if $a_i>0$} \\
              -   &\mbox{if $b_i>0$} 
            \end{array}
    \right.  
\mbox{ for every } i\in E'
\]
We apply Proposition~\ref{contractible} with this $R'$ to 
%the affine oriented matroid
$\,(E \backslash E'', \L / E'', g)$.
Then $B'$ is the face poset of $(B_\M)_{\leq \x^{\bf a}\y^{\bf b}}$,
which is therefore contractible.
\endproof

{\smallskip\noindent}
The oriented matroid ideal $O$ is square-free, and hence
is the Stanley-Reisner ideal 
of a simplicial complex $\Delta_{\M}$
on $2n$ vertices $\{1,\ldots,n,1', \ldots,n'\}$,
whose faces correspond to square-free
monomials of $\k[\x,\y]$ that  do not belong to $O$, i.e.,
$$ \{i_1, \ldots, i_k, j_1', \ldots, j_m'\}\in \Delta_{\M}
\quad  \mbox{ if and only if } \quad
x_{i_1}\ldots x_{i_k} y_{j_1} \ldots y_{j_m} \notin O.$$
In what follows we give a geometric description
of that simplicial complex.

\begin{lemma}       \label{i,i'}
We have $ F\cap\{i, i'\}\neq \emptyset$ for 
any facet $F$ of $\Delta_{\M}$  and
$\, i\in\{1,\ldots,n\}$.
\end{lemma}

\proof Let $F$ be a face of $\Delta_{\M}$ such that 
$F\cap\{i, i'\}=\emptyset$. Suppose that neither $F'=F\cup \{i\}$
nor $F''=F\cup \{i'\}$ is a face of $\Delta_{\M}$.
Then there exist cocircuits $Z', Z'' \in B_\M$ such that
\[
 Z'_i = +, \;\;\; (Z')^+ \setminus\{i\} \subseteq 
                   \{1\leq j \leq n : j\in F\} \cup \{ g\}, \;\;\;
       (Z')^- \subseteq \{1\leq j \leq n : j'\in F\}   
\]\[
Z''_i = -, \;\;\;(Z'')^+ \subseteq \{1\leq j \leq n : j\in F\} \cup \{g\}, \;\;
    (Z'')^- \setminus\{i\} \subseteq \{1\leq j \leq n : j'\in F\}.   
\]
By the strong elimination axiom applied to $(Z', Z'', i, g )$,
there is a cocircuit $Z$ such that 
$
Z_i = 0, \; Z_g=+, \; Z^+ \subseteq (Z')^+ \cup (Z'')^+, \;
   Z^- \subseteq (Z')^- \cup (Z'')^-.
$
Thus $Z \in B_\M$, and the monomial $m_{xy}(F)$ 
is divisible by $m_{xy}(Z) \in O$. This
contradicts  $F \in \Delta_{\M}$.
\endproof 

Suppose now that an affine oriented matroid $\M=(E, \L, g)$ 
is a single element extension of  the matroid 
$\M\backslash g=(E\backslash g, \L \backslash g)$
by an element $g$ in {\it general position,}
in the sense of  \cite[Proposition~7.2.2]{BLSWZ}. For
the affine arrangement $\A$ in Section 1 or Example 2.1, 
this means that $\A$ has no vertices at infinity.
In such a case, Theorem~\ref{res1} implies the following properties of $O$.
In the rest of this section we denote by  $r$ the rank 
of~$\M\backslash g$.

\begin{corollary}    \label{CM}
 $\field[\Delta_{\M}]=\k[\x,\y]/O$ 
is a Cohen-Macaulay ring of dimension $2n-r$. 
\end{corollary}
\proof
Since $\rank(\M \backslash g)=r$,
every $(n-r+1)$-element subset $\{i_1,\ldots,i_{n-r+1} \}$
 of $\{1, \ldots, n\}$ contains the support
of a (signed) cocircuit. This implies that every monomial
of the form $x_{i_1} \ldots x_{i_{n-r+1}}y_{i_1} \ldots y_{i_{n-r+1}}$ 
belongs to $O$. The variety defined by these monomials is
a subspace arrangement of codimension $r$. Hence $O$
has codimension $\leq r$, which means that the ring
$\, \field[\Delta_{\M}] = \k[\x,\y]/O\,$ has
Krull dimension $\,\leq 2n-r$.
By Theorem~\ref{res1}, the bounded complex $B_\M$
supports a minimal free resolution of $O$, and therefore
$$\depth(\field[\Delta_{\M}]) \,\, =\,\,
2n-\mbox{(the length of this resolution)} \,\, = \,\, 2n -r. $$
Hence $\depth(\field[\Delta_{\M}])=
\dim(\field[\Delta_{\M}])=2n-r$,
and $\field[\Delta_{\M}]$ is Cohen-Macaulay. \endproof

\begin{corollary}             \label{reg_seq}
$\{ x_1-y_1, \ldots, x_n-y_n \}$ is a regular sequence on 
$\field[\Delta_{\M}]=\k[\x,\y]/O$.
\end{corollary}

\proof
Since $\field[\Delta_{\M}]$ is
Cohen-Macaulay,
it suffices to show that $\{ x_1-y_1, \ldots, x_n-y_n \}$ is a part of
a linear system of parameters (l.s.o.p.).
This follows from
Lemma~\ref{i,i'} and the l.s.o.p.~criterion
due to Kind and Kleinschmidt \cite[Lemma III.2.4]{St96}. \endproof

Consider any signed circuit $C = (C^+,C^-)$ of our oriented matroid 
such that $g$ lies in $C^{-}$.  By 
the {\it general position\/} assumption on $g$, the complement
of $g$ in that circuit is a basis of the underlying matroid.
We write $P_C$ for the
ideal generated by the variables $x_i$ for each
$i \in C^+$  and the variables
$y_j$ for each  $j \in C^-\backslash  \{ g\}$.

\begin{proposition} 
\label{primeofO}
The minimal prime decomposition of the oriented matroid ideal equals
$\, O \,=\, \bigcap_C P_C \,$ where the intersection is over all
circuits $C$ such that  $g \in C^- $.
\end{proposition} 
\proof
The right hand side is easily seen to contain the left hand side.
For the converse it suffices to divide by the regular
sequence $x_1- y_1,\ldots,x_n-y_n$ and note that the
resulting decomposition for the matroid ideal $M$ is easy
(Remark \ref{easyremark}).
\endproof

Our final result relates the ideal $O$ to
matroid polytopes and their triangulations. The monograph of 
Santos \cite{San} provides an excellent state-of-the art introduction.
We refer in particular to \cite[Section 4]{San}, where
Santos introduces triangulations of Lawrence (matroid) polytopes,
and he shows that these are in bijection with one-element liftings
of the underlying matroid. Under matroid duality, one-element liftings
correspond to one-element extensions. In our context, these extensions
correspond to adding the special element $g$ which plays
the role of the pseudohyperplane at infinity.  From Santos' result we 
infer the following theorem.

\begin{theorem}
\label{lawlaw}
The oriented matroid ideal $O$ is the Stanley-Reisner ideal 
of the triangulation of the Lawrence matroid polytope
induced by the lifting dual to the extension by $g$. In particular,
$O$ is the Stanley-Reisner ideal of a triangulated ball.
\end{theorem}

The second assertion holds because
lifting triangulations of matroid polytopes are triangulated balls
and, by Santos' work, every triangulation of
a Lawrence matroid  polytope is a lifting triangulation.
We remark that it is unknown whether arbitrary triangulations
of  matroid polytopes are topological balls \cite[page 7]{San}.

\bigskip

\section{Matroid Ideals} 
\label{matroid_section}

Let $\matroid$ be an (unoriented) matroid on the set 
$\{1, \ldots, n\}$ and let $L$ be its lattice of 
flats. We encode $\matroid$ by the
{\it matroid ideal} $M$  generated by
the monomials $m_x(F) = \prod_{i: i\notin F} x_i$
for every proper flat $F \in L$. The minimal generators
of $M$ are the squarefree monomials representing 
cocircuits of $\matroid$, that is, the monomials $m_x(H)$ 
where $H$ runs over all hyperplanes of $\matroid$.
Equivalently, $M$ is the Stanley-Reisner ideal of the simplicial
complex of independent sets of the dual matroid $\matroid^*$.
This explains what happens when we substitute
$\,y_i \mapsto x_i \,$ in Proposition~\ref{primeofO}:

\begin{remark}
\label{easyremark}
The matroid ideal $M$ has the minimal prime decomposition
$$ M \quad = \quad \bigcap_{B \textrm{ \rm basis of } \matroid}
\langle \, x_i \,\,| \,\, i \in B \,\rangle. $$
\end{remark}

The following characterization of our ideals
can serve as a definition of the word ``matroid''.
It is a translation of the (co)circuit axiom into commutative algebra.

\begin{remark} \rm
A proper square-free monomial
ideal $M$ of $\k[\x]$ is a {\it matroid ideal} if and only if
for every pair of monomials $m_1, m_2 \in M$ and any $i\in\{1, \ldots,n\}$
such that $x_i$ divides both $m_1$ and $m_2$, the monomial
$\lcm(m_1, m_2) / x_i $ is in $M$ as well.
\end{remark}

\vskip .1cm

Matroid ideals have been studied since the earliest days of
combinatorial commutative algebra, as a paradigm
for shellability and Cohen-Macaulayness. Stanley computed their
Betti numbers  in \cite[Theorem~9]{St77}. The purpose
of this section is to construct an explicit minimal 
$\field[\x]$-free resolution for any matroid ideal $M$.

We note that Reiner and Welker \cite{RW} used the term
``matroid ideal'' for the square-free monomial ideals which are 
Alexander dual to our ideals. The matroid ideals in \cite{RW} have 
linear resolution but are generally not Cohen-Macaulay, while our
matroid ideals are Cohen-Macaulay but their resolution
(given below) is generally not linear. In particular,
the Alexander dual of a matroid ideal is completely different
from the matroid ideal of the dual matroid.

We first consider the case where $\matroid$ is an {\it orientable matroid}.
This means that there exists an  oriented matroid $\M$
whose underlying matroid is $\underline\M$.
Let  $\L$ be the set of covectors of a single element extension of
 $\M$ by an element $g$ in general position 
(see \cite[Proposition~7.2.2]{BLSWZ}). 
%Then $\widetilde{\M}$
%is an oriented matroid on the set $E=\{1,\ldots,n\}\cup\{g\}$, and, 
Consider the affine oriented matroid $\widetilde{\M}=(E, \L, g)$, where 
$E=\{1,\ldots,n\}\cup\{g\}$,
and its bounded complex $B_{\widetilde{\M}}$.
Note that, for each sign vector $Z$ in $B_{\widetilde{\M}}$,
the zero set $Z^0$ is a flat in $L$. Moreover, by the genericity
hypothesis on $g$, all flats arise in this way.
We label each cell $Z $ of the bounded complex
$B_{\widetilde{\M}}$ by the monomial
$\,m_x(Z)=\prod \{ x_i \,: \,1\leq i\leq n \mbox{ and }  Z_i \neq 0 \}$.

\begin{theorem}  \label{res2} 
Let $M$ be the matroid ideal of an orientable matroid. Then 
the bounded complex $B_{\widetilde{\M}}$ 
of any corresponding affine oriented matroid 
is an $M$-complex, and its cellular complex
$C_\bullet(B_{\widetilde{\M}}, M)$ gives a minimal 
free resolution of~$M$ over $\k[\x]$.
\end{theorem}

\proof
Let  ${\bf a} = (a_1,\ldots,a_n) \in\N^{n}$ and  consider $M$-essential
subcomplex $(B_{\widetilde{\M}})_{\leq \x^{\bf a}}$. This complex 
(if not empty) is the bounded complex of the contraction of
$(E,\L, g)$  by $\,\{1\leq i\leq n: a_i= 0\}$, and, 
hence is acyclic by Proposition~\ref{contractible}. Since
$m_x(Z')$ is a proper divisor of $m_x(Z)$ whenever
$Z' < Z$, and $Z', Z \in B_{\widetilde{\M}}$, it follows that
$B_{\widetilde{\M}}$ is an $M$-complex.
\endproof

We remark that  $C_\bullet(B_{\widetilde{\M}}, M)$
is obtained from the complex $C_\bullet(B_{\widetilde{\M}}, O)$,
where $O$ is the oriented matroid ideal of
$\widetilde{\M}=(E,\L,g)$,
by specializing $y_i$ to $x_i$ for all $i$.
Hence Theorem~\ref{res1} and Corollary \ref{reg_seq} 
give a second proof of Theorem~\ref{res2}.

\begin{corollary}
\label{happiness}
The $\N^n$-graded  Hilbert series of any matroid ideal $M$ equals
\begin{equation}
\label{matroidhilbert}
 \biggl(
 \sum_{F\in L} \mu_L (F, \hat{1})
    \cdot \prod \{x_j : j\notin F \} \biggr)
/ \prod_{i=1}^n (1-x_i)
\end{equation}
where $L$ is the lattice of flats of $\matroid$, and $\mu_L$
is its M\"obius function.
\end{corollary}

There are several ways of deriving this corollary.
First, it follows from \cite[Theorem~9]{St77}.
A second possibility is to observe that
the geometric lattice $L$ coincides with the lcm lattice
(in the sense of \cite{GPW}) of the ideal $M$,
and then \cite[Theorem~2.1]{GPW} implies the claim.
Finally, in the orientable case, Corollary  \ref{happiness}
follows from Theorem~\ref{res2} and the oriented matroid version of
Zaslavksy's face-count formula.

\begin{proposition} 
{\rm (Zaslavsky's Formula) \cite{Zas}, \cite[Theorem~4.6.5]{BLSWZ}} \
\label{Mobius}
The number of bounded regions of a rank $r$ affine oriented matroid 
${\widetilde{\M}}\!=\!(\!E,\L,g\!)$
equals~$(\!-1\!)^{r}\!\mu_{L}(\hat{0},\! \hat{1})$.
\end{proposition}  

\smallskip

We next treat the case of non-orientable matroids. It would be
desirable to construct an $M$-complex for an arbitrary matroid
ideal $M$, and to explore the ``space'' of all possible
$M$-complexes. Currently we do not know how to construct them.
Therefore we introduce a different technique for
resolving $M$ minimally.

Let $P$ be any graded poset which has a unique minimal element 
$\hat{0}$ and a unique maximal element $\hat{1}$.
(Later on, we will take $P$ to be the order dual of our geometric lattice $L$).
Let $\O(P)$ denote the order complex of $P$, 
that is, the simplicial complex whose
simplices  $[F_0, F_1, \ldots, F_i]$ are decreasing chains
$\, \hat{1} > F_0>F_1 > \ldots > F_i > \hat{0}$.
For  $F \in P$ denote by $\O(F)$  the order complex of the lower interval
$[\hat{0},F]$. Note that 
$\dim\O(F) = \rank(F)-2$.
Let  $C_i(\O(F))$ be the $\k$-vector space of $i$-dimensional
chains of $\O(F)$,
and let 
$$ 0 \longrightarrow
 C_{\rank(F)-2}(\O(F)){\longrightarrow}\, \ldots
%\,{\overset{\partial_4}\longrightarrow}\, \C_3(\Delta,\M)
\,{\overset{\partial_2}\longrightarrow}\, C_1(\O(F))
\,{\overset{\partial_1}\longrightarrow}\, C_0(\O(F))
\,{\overset{\partial_0}\longrightarrow}\, C_{-1}(\O(F))
\,{\longrightarrow}\, 0\,.
$$
be the usual (augmented) chain complex, i.e., the differential is given by 
$$\partial_i[F_0, F_1, \ldots, F_i] = \sum_{j=0}^i (-1)^j 
      [F_0, \ldots, \widehat{F_j}, \ldots, F_i] \,\,\,\,
\mbox{ for } i>0 \;\;
\mbox{ and }
\partial_0[F_0]=\emptyset.
$$
Denote by $Z_i(\O(F))=\ker (\partial_i)$ 
the space of  $i$-cycles, and by $\widetilde{H}_i(\O(F))$
the $i$th (reduced) homology of $\O(F)$.
For relevant background on poset homology see \cite{Bj}.

For each pair $F, F' \in P$  such that $\rank(F)-\rank(F')=1$, 
we define a map
$$\phi: C_i (\Delta_{F}) \longrightarrow C_{i-1}(\Delta_{F'}) \mbox{ by }
[F_0, F_1, \ldots,  F_i] \mapsto
\left\{ \begin{array}{ll}
              0,   &\mbox{if $F_0 \neq F'$} \\
              {[}F_1, \ldots, F_i{]}    &\mbox{if $F_0=F'$.} \\
            \end{array}
    \right.  
$$
The map $\phi$ is zero unless $F'<\!\!\! \cdot \, F$
(in words: $F$  covers $F'$).
Note that $\partial\circ \phi = -\phi\circ \partial$, 
and hence the restriction of $\phi$ to cycles gives a map 
$\phi: Z_i(\O(F)) \longrightarrow Z_{i-1}(\O(F'))$.
Combining  these maps together we obtain a complex of $\k$-vector spaces:
\begin{eqnarray*}  
 \ZZ(P):  && 0  \, \longrightarrow \,
 Z_{r-2}(\O(P))  \,\,\,{\overset{\phi} \longrightarrow}
 \bigoplus_{\rank(F)=r-1} \!\!\!\! Z_{r-3}(\O(F))
\,\,\, {\overset{\phi}
\longrightarrow}\, \,\,
\cdots  \\
&& \cdots \,\,\,{\overset{\phi}\longrightarrow}\, 
\bigoplus_{\rank(F)=2} \!\!\!  Z_0(\O(F)) \,\,\,
{\overset{\phi}\longrightarrow}\, 
\bigoplus_{\rank(F)=1} \!\!\! Z_{-1}(\O(F)) \,\,
\,{\longrightarrow}\, \,\, \k\,.    \nonumber
\end{eqnarray*}
The complex property $\phi^2 = 0$ is verified by direct calculation
using the equation (\ref{sum_i}) stated below.
Let $P_{(j)}$ denote the poset obtained from $P$ by
removing all rank levels $\geq j$, and let $\O(P_{(j)})$ be the order complex
of $P_{(j)}\cup\{\hat{1}\}$.

\begin{proposition}   
\label{prop:k-exact}
The complex $\ZZ(P)$ is exact if
$\,\widetilde{H}_{i} (\O(P_{(i+3)}))= 0 \,$ 
for all $ i\leq r-3$.
\end{proposition}

To prove Proposition~\ref{prop:k-exact} we need some notation.
If $x\in\bigoplus_{\rank(F)=i}Z_{i-2}(\O(F))$, we denote its $F$-component
by $x_F$. For a simplex $\sigma=[F_0,F_1,\ldots,F_i]$ we also 
write  $\sigma=F_0 * [F_1, \ldots, F_i]$, and 
the operation ``$*$'' extends to $\k$-linear combinations.

\begin{remark}  \label{cycle}
{\rm
Suppose that $z\in C_{i}(\O(P_{(i+2)}))$.
Then $z$ can be expressed as 
$$
z \quad =\sum_{\rank(F')=i+1} \!\! F'* y_{F'}  \quad
=\, \sum_{{\rank(F')=i+1}}\sum_{ F''<\!\!\cdot F'} \! F'* F'' * x_{F', F''},
$$
where $y_{F'} \in C_{i-1}(\O(F'))$ and  $x_{F', F''} \in C_{i-2}(\O(F''))$.
Its boundary equals
\begin{eqnarray*}
\partial(z) &=  \,\,\,\,
\sum_{\rank(F'')=i}  F'' * \sum_{F'\cdot\!\!> F''}x_{F', F''} & \\
&- \,\sum_{\rank(F')=i+1} F' * \sum_{F''<\!\!\cdot F'}x_{F', F''}& +
\sum_{F', F''}
      F'  *   F'' * \partial(x_{F', F''}).
\end{eqnarray*}
We conclude that $z$ is a cycle if and only if the following
conditions are satisfied:
\begin{eqnarray}
\sum_{F'\cdot\!\!> F''}x_{F', F''}=0 && 
\hbox{for all} \,\, F'' \,\, \hbox{with} \,\, \rank(F'')=i 
                      \label{sum_i}\\
\sum_{F''<\!\!\cdot F'}x_{F', F''}=0 && 
\hbox{for all} \,\, F' \,\, \hbox{with} \,\, \rank(F')=i+1
                      \label{sum_i+1}\\
\partial(x_{F', F''})=0 && \hbox{for all} \,\,\, F', F'' \mbox{ such that } 
                      F''<\!\!\!\cdot F'.   \label{boundary}
\end{eqnarray} 
}
\end{remark}

{\smallskip\noindent {\it Proof of Proposition~\ref{prop:k-exact}: }}
To show that $\ZZ(P)$ is exact, consider 
$y=(y_{F'})\in\bigoplus_{\rank(F')=i+1} Z_{i-1}(\O(F'))$ such that 
$\phi(y)=0$. There are several cases:
If $i=r-1$, then $y=y_{\hat{1}}$ can be expressed as
$\sum_{\rank(F)=r-2} F * x_F$, where $x_F \in C_{r-3}(\O(F))$.
Then  $0=\phi(y)_F =x_F$ and therefore $y=0$.
Hence the leftmost map $\phi$ is an inclusion.  

Let  $0 < i<r-1$ and define
$\,z \, =\,
\sum_{\rank(F')=i+1} F'* y_{F'} \in C_{i}(\O(P_{(i+2)}))$.
We claim that $z$ is a cycle, that is,
$z\in Z_{i}(\O(P_{(i+2)}))$.
Indeed, if $i>0$, then 
 $y_{F'}$ can be expressed as
$\sum_{ F''<\!\!\cdot F'} F'' * x_{F', F''}$, where
$x_{F', F''} \in C_{i-2}(\O(F''))$. Hence
\begin{eqnarray*}
(\phi(y))_{F''} = \sum_{F'\cdot\!\!> F''}x_{F', F''} &
  \forall F'' \mbox{ with }
\rank(F'')=i,  \mbox{ and }& \\
\partial(y_{F'})= \sum_{ F''<\!\!\cdot F'} x_{F', F''}  \, -  \!
                  \sum_{ F''<\!\!\cdot F'} F'' * \partial (x_{F', F''}) 
& \forall F' \mbox{ with }  \rank(F')=i+1. &
\end{eqnarray*}
Since $\phi(y)=0$ and $\partial(y_{F'})=0$ for any $F'$ of rank $i+1$,
we infer that $z$ satisfies conditions (\ref{sum_i})--(\ref{boundary})
in  Remark \ref{cycle}, and therefore  is a cycle.
In the case $i=0$ the proof is very similar.
 Now, if $i=r-2$ then 
$z\in Z_{r-2}(\O(P))$, and  
  $\phi(z)=\phi(\sum F'* y_{F'})=(y_{F'})=y$.
Hence we are done in this case.
If $i<r-2$, then, since
$Z_{i}(\O(P_{(i+2)}))\subseteq  Z_{i}(\O(P_{(i+3)}))$,
and  $\widetilde{H}_{i} (\O(P_{(i+3)}))= 0$, it follows that 
 there exists $w \in  C_{i+1}(\O(P_{(i+3)}))$ such that 
$\partial(w)=z$. Express $w$ as $\sum_{\rank(F)=i+2} F* v_F$, where
$v_F \in C_i(\O(F))$. 
Since 
$z=\partial(w)=\sum_{\rank(F)=i+2} v_F + 
\sum_{\rank(F)=i+2} F*\partial(v_F)$, we conclude that
$\partial(v_F)=0$ for all $F$ of rank $i+2$, and that 
$\sum_F v_F = z = \sum_{F'} F' * y_{F'}$. Thus
$v=(v_F) \in \bigoplus_{\rank(F)=i+2} Z_{i}(\O(F))$, and 
$\phi(v)=y$.
\endproof

\begin{corollary}
\label{cmposet}
If $P$ is a Cohen-Macaulay poset, then $\ZZ(P)$ is exact.
\end{corollary}
\proof If $\O(P)$ is Cohen-Macaulay,
then  $\O(P_{(i)})$ is Cohen-Macaulay for every $i$ (see
 \cite[Theorem~4.3]{St79}). This means that
all homologies of $\O(P_{(i)})$ vanish, except possibly for the top
one. Thus the conditions of Proposition~\ref{prop:k-exact} are satisfied.
\endproof

Suppose now that every atom $A$ of $P$ is labeled by a monomial $m_A\in\k[\x]$.
The {\it poset ideal\/} $I_P$ is the ideal generated by these monomials.
Associate with every element $F$ of $P$ a monomial $m_F$
as follows:
\[
m_F:=\lcm\{m_A : \rank(A)=1,  A\leq F\} \;\: \mbox{ if }
                 F\neq \hat{0}\; \; \mbox{ and } m_{\hat{0}}:=1.
\]
We say that the labeled poset  $P$ is {\it complete\/} if all monomials
$m_F$ are distinct, and  for every
${\mathbf a}\in \N^{n}$ the set $\{F\in P : \deg(m_F)\leq {\mathbf a}\}$
has a unique maximal element.

We identify the principal ideal $\langle m_F \rangle$ with the free 
$\N^n$-graded $\k[\x]$-module of rank 1 with generator in degree $\deg m_F$.
If $F, G \in P$ and $F<  G$, then $m_F$ is a divisor of $m_G$.
Thus there is an inclusion of the  corresponding ideals
$\, \langle m_G \rangle \longrightarrow \langle m_F \rangle $.
Recall that there is a complex $\ZZ(P)$ of $\k$-vector spaces 
associated with $P$. 
Tensoring summands of this complex with the ideals 
$\{\langle m_F \rangle: F\in P\}$,
we obtain a complex of $\N^n$-graded free $\k[\x]$-modules :
\begin{equation}      \label{C-sequence}
\C(P)\,\, = \,\, \bigoplus_{F\in P} Z_{\rank(F)-2}(\O(F))\otimes_{\k} 
   \langle m_F \rangle
\quad \mbox{ with differential } \partial=\phi\otimes i.
\end{equation}

\begin{theorem}   \label{S-exact}
Suppose that the labeled poset $P$ is complete and that the homology
$\widetilde{H}_i(\O(F_{(i+3)}))$  vanishes for any $0\leq i \leq r-3$  and 
any $F\in P$ of rank $\geq i+3$,  then $(\C(P), \partial)$
is a minimal $\N^n$-graded free $\k[\x]$-resolution of the poset 
ideal $I_{P}$.
\end{theorem}
\proof
$(\C(P), \partial)$ is a complex of $\N^n$-graded
free $\k[\x]$-modules. To show that it is a resolution, we have to
check that for any ${\mathbf a} \in \N^n$, the ${\mathbf a}$th graded
component  $\,(\C(P), \partial)_{\mathbf a} \,$
is an exact complex of $\k$-vector spaces. 
Let ${\mathbf a} \in \N^n$, and  let $F\in P$ be the maximal
element among all elements $G\in P$ such that $\deg(m_G)\leq{\mathbf a}$.
Such an element $F$ exists since the labeled poset $P$ is complete.
Then $(\C(P), \partial)_{\mathbf a}$ is isomorphic to 
the complex $\ZZ([\hat{0}, F])$ of the poset $[\hat{0}, F]$, 
and hence is exact over $\k$ (by Proposition~\ref{prop:k-exact}).
Thus $(\C(P), \partial)$ is exact over $\k[\x]$.
Finally, since $m_F$ and $m_G$ are distinct monomials for any pair
$F<\!\!\!\cdot G$, the resolution $(\C(P), \partial)$ is minimal.
\endproof

{}From Corollary \ref{cmposet} we obtain:

\begin{corollary}  \label{CM-exact}
If $P$ is a complete labeled poset such that
 every lower interval of $P$ is Cohen-Macaulay, then 
$(\C(P), \partial)$
is a minimal $\N^n$-graded free resolution of $I_{P}$.
\end{corollary}

Returning to our matroid $\matroid$,
let $P$ be a  lattice of flats ordered by 
 reverse inclusion.  Hence $P$ is the order dual of the geometric lattice 
$L$ above. In particular,  $\hat{0}$ corresponds
to the set $\{1,2, \ldots, n\}$, and $\hat{1}$
 corresponds to the empty set. Label each atom $H$ of $P$
(that is, hyperplane of $\matroid$) by the monomial $m_x(H)$ as in the 
beginning of this section. Identifying the variables $x_i$ with
the coatoms of $P$, we see that $m_x(H)$ is the product
over all coatoms not above $H$. Then $P$ is a complete labeled poset
and its poset ideal $I_P$  is precisely the matroid ideal $M$.
Moreover, all lower intervals of the poset $P$
are Cohen-Macaulay \cite[Section 8]{St77}.
{}From  Corollary \ref{CM-exact} we obtain the
following alternative to Theorem~\ref{res2}. 

\begin{theorem}
Let $\matroid$ be any matroid. Then  the complex
$(\C(P), \partial)$
is a minimal $\N^n$-graded free $\k[\x]$-resolution of the  
matroid ideal~$M$.
\end{theorem}

The two resolutions presented in this section provide a  syzygetic realization
of Stanley's formula \cite[Theorem~9]{St77} for the Betti numbers
of matroid ideals. That formula states that the number of minimal $i$-th 
syzygies of $\k[\x]/M$ is equal to
$$  \beta_i(M) 
\quad = \quad \sum_{F} |\,\mu_L (F, \hat{1})\,|, $$
where the sum is over all flats $F$ of 
corank $i$ in $\matroid$. The generating function 
\begin{equation}
\label{stanley}
\psi_{\matroid}(q)=\sum_{i=0}^{\rank(M)}  \beta_i(M) \cdot q^i  \quad = \quad
\sum_{F \textrm{ flat of} \,\matroid}
|\,\mu_L (F, \hat{1})\,| \cdot q^{{\rm corank}(F)} 
\end{equation}
for the Betti numbers of $M$ 
is called the {\it cocharacteristic polynomial} of $\matroid$.
In the next two sections we will examine this polynomial
for some special matroids.

\bigskip

\section{Unimodular toric arrangements}

A {\it toric arrangement} is a hyperplane arrangement which lives
on a torus  $\T^d$ rather than in $\R^d$. One construction of
such arrangements appears in recent work
of Bayer, Popescu and Sturmfels \cite{BPS}. 
Experts on geometric combinatorics might appreciate
the following description:
Fix a unimodular matroid $\matroid$, form the associated
tiling of Euclidean space by zonotopes 
\cite[Proposition~3.3.4]{Wh}, dualize to get an 
infinite arrangement of hyperplanes, and 
divide out by the group of lattice translations.

Here is the same construction again, but now in slow motion. Fix a central
hyperplane arrangement $\C=\{H_1, \ldots, H_n\}$ in
$\R^d$ where $H_i=\{ v \in \R^d : h_i \cdot v=0\}$ for some
$h_i\in \Z^d$. Let $L$ denote the intersection lattice of $\C$
ordered by reverse inclusion.
 We assume that $\C$  is {\it unimodular,}
which means that the $d\times n$ matrix $(h_1, \ldots, h_n)$
has rank $d$, and all its $d \times d$-minors lie in the set
$\{0, 1, -1\}$. We retain this hypothesis throughout this section.
 See \cite{Wh} and \cite[Theorem~1.2]{BPS} for details on
unimodularity. 
The set of all integral translates of hyperplanes of~$\C$,
$$ H_{ij}=\{v\in \R^d :  h_i\cdot v= j\} \quad \mbox{ for } 
i\in\{1, \ldots, n\}  \mbox{ and } j\in \Z, $$
forms an infinite arrangement $\widetilde{\C}$ in $\R^d$.
The unimodularity hypothesis is equivalent to saying
that the set of vertices of $\widetilde{\C}$
is precisely the lattice  $\Z^d$, that is, no
new vertices can be formed by intersecting the
hyperplanes $H_{ij}$.   Define the 
{\it unimodular toric arrangement\/} $\widetilde{\C}/ \Z^d$ 
to be the set of images of the $H_{ij}$ in the torus $\T^d=\R^d/\Z^d$.

Slightly abusing notation, we refer to these images 
as hyperplanes on the torus. The images of cells
of~$\widetilde{\C}$ in $\T^d$ are called {\it cells} of 
$\widetilde{\C}/\Z^d$. These cells form a
cellular decomposition of $\T^d$.  
 Denote by  $f_i=f_i(\widetilde{\C}/\Z^d)$ 
the number of $i$-dimensional cells in this decomposition. 
The next result concerns the {\it $f$-vector\/}
$(f_0,f_1,\dots,f_d)$ of $\widetilde{\C}/\Z^d$.

\begin{theorem}     \label{f-vector}
If $\,\widetilde{\C}/\Z^d \,$ is a unimodular toric arrangement,  then 
$$\sum_{i=0}^d f_i(\widetilde{\C}/\Z^d) \cdot q^i
\,\,\, = \,\,\, \psi_\C(q), \quad \mbox{ where } \,\,
\psi_\C(q) \, = \, \sum_{F\in L}\mu_L(F, \hat{1}) \cdot (-q)^{\dim F} $$
is the  cocharacteristic polynomial of the underlying
hyperplane arrangement $\,\C$.
\end{theorem}

\proof
Choose a vector $w\in \R^d$ which
  is not perpendicular to any 1-dimensional cell of the arrangement $\C$.
Consider the  affine hyperplane $\{ \,v\in \R^d \,:\, w \cdot v=1\}$.
Let $\A=\C\cap H$ be a restriction of  $\C$ to $H$. 
Then $\A$ is an affine arrangement in $H$. For any  $i\geq 0$,
there  is a one-to-one correspondence between the $(i-1)$-dimensional 
bounded cells of $\A$ and the $i$-dimensional cells of toric arrangement 
$\widetilde{\C}/\Z^d$. To see this, consider the  cells in the
infinite arrangement $\widetilde{\C}$ whose minimum with respect to the 
linear functional $\, v \mapsto \,w \cdot v \,$ is attained at the origin. 
These cells form a system of representatives modulo the $\Z^d$-action. But they
are also in bijection with the bounded cells of $\A$.
Using Proposition~\ref{Mobius} (see also Example \ref{central}), we conclude
$$
f_i(\widetilde{\C}/\Z^d)
\quad = \quad  f_{i-1}(B_\A) \quad = \quad 
(-1)^i \cdot \!\! \sum_{\dim(F)=i}\mu_L(F, \hat{1}), $$  
where the sum is over elements of $L$ of corank $i$. 
This completes the proof.
\endproof

Theorem~\ref{f-vector} was found independently by Vic Reiner
who suggested that we include the following alternative proof.
His proof has the advantage that it does not rely
on Zaslavsky's formula. 

\medskip

\noindent {\sl Second proof of Theorem~\ref{f-vector}: }
Starting with the unimodular toric arrangement $\torus$,
for each intersection subspace $F$ in the intersection
lattice $L$, let $T_F$ denote the subtorus obtained by
restricting $\torus$ to $F$.  So $T_0$ is just $\torus$
itself, and $T_1$ is not actually a torus but rather a point.
Our assertion is equivalent to
\begin{equation}
\label{vicsequation}
\mu(F,1) \quad = \quad  (-1)^{\dim F} \cdot \#\{\text{max cells in }T_F\}
\end{equation}
Let  $\mu'(F)$ denote the right-hand side above.
By the definition of the M\"obius function of a poset,
the equation (\ref{vicsequation}) is equivalent to
$$
\sum_{F \leq G \leq 1} \mu'(G) \,\,\, = \,\,\, \delta_{F,1}.
\qquad \qquad \hbox{(Kronecker delta)}
$$
The left hand side of this equation is the
(non-reduced) Euler characteristic of $T_F$.
This is $0$ since $T_F$ is a torus, unless
$F=1$ so that $T_F$ is a point and then it is $1$.
\endproof

\smallskip

We remark that Theorem~\ref{f-vector} can be generalized to
arbitrary toric arrangements $\widetilde{\C}/\Z^d$
without the unimodularity hypothesis. The face count formula
is a sum of local M\"obius function values
over all (now more than one) vertices of $\widetilde{\C}/\Z^d$.
That generalization has interesting applications to hypergeometric functions,
which  will be the subject of a future publication. 
The enumerative applications 
in the next section all involve unimodular arrangements, so we restrict 
ourselves to this case. We shall need the following recursion
for computing cocharacteristic polynomials.

\begin{proposition} \label{prop:psi-recurrence}
Let $H$ be a  hyperplane of the arrangement~$\C$.  Then
$$
\psi_\C (q) \quad = \quad \psi_{\C\cap H}(q) \, + \,q 
\cdot \sum_c \psi_{\C/c}(q),
$$
where the sum is over all lines $c$ of the arrangement $\C$
that are not contained in~$H$.
\end{proposition}

The lines $c$ of the arrangement $\C$ are
the coatoms of the intersection lattice $L$. The arrangement
$\C / c$ is the hyperplane arrangement
$\{ \, H_i /  c \,\, : \,\, c \in H_i \, \} \,$ in
the $(d-1)$-dimensional vectorspace $\R^d/c $.
Note that if $c$ is a simple intersection, that is,
if  $c$ lies on only $d-1$ hyperplanes $H_i$, then
$\,\psi_{\C/c}(q) \, = (1+q)^{d-1} $. 
Note that Proposition~\ref{prop:psi-recurrence}, 
together with the condition $\psi_\C(q)=1$ for
the $0$-dimensional arrangement~$\C$,
uniquely defines the cocharacteristic polynomial.

\proof 
The intersection lattice $L$ of any central hyperplane arrangement 
$\C$ is  {\it semi-modular}, that is, if  both $F$ and $G$ cover 
$F \wedge G $,
then $F \vee G $ covers both $F$ and $G$, 
see~\cite[Section~3.3.2]{Stanley:EC}.
The assertion follows from the 
relation~\cite[formula~3.10.(27)]{Stanley:EC}
for the M\"obius functions of any semi-modular lattice.
\endproof

In the remainder of this section we review the
algebraic context in which unimodular toric arrangements
arise in \cite{BPS}. This provides a
Gr\"obner basis interpretation for our proof of Theorem \ref{f-vector}
and it motivates  our enumerative results in Section~5.

Denote by $B$ the $n\times d$ matrix whose rows are $h_1, \ldots, h_n$.
All $d \times d$-minors of $B$ are $-1,0$, or $+1$.
The {\it unimodular Lawrence ideal} of $\C$ is the binomial prime ideal
$$
J_\C \,\,\, := \,\,\,
\langle \,\x^a \y^b -\y^a \x^b \,\, | \,\,
a, b \in \N^n \, , \,\,
 a-b \in \mathrm{Image} (B) \,\rangle \quad \hbox{in} \quad
\k[x_1, \ldots, x_n, y_1, \ldots, y_n]. $$
The main result of \cite{BPS} states that the toric arrangement
$\,\C / \Z^d \,$ supports a cellular resolution of $J_\C$.
In particular, the Betti numbers of the unimodular Lawrence ideal
$J_\C$ are precisely the coefficients of the cocharacteristic polynomial
$\,\psi_{\C}(q) $.

The construction in the proof of Theorem~\ref{f-vector}  has a 
Gr\"obner basis interpretation. Indeed, the generic vector $w \in \R^d$
defines a term order $\succ$ for the ideal $J_\C$ as follows:
$$ 
\x^a \y^b \,\succ \, \y^a \x^b 
\quad \hbox{if} \quad
 \,\,\
a-b\,=\,B \cdot u \,\,\, \hbox{for some $u \in \R^d$
with}\,\,\,\,w\cdot u > 0. $$
It is shown in \cite[\S 4]{BPS} that
the initial  monomial ideal $in_\succ(J_\C)$ of $J_\C$ with respect to 
these weights
is the oriented matroid ideal associated with the restriction of
the central arrangement $\C$ to the affine hyperplane 
$\{ \,v\in \R^d \,:\, w \cdot v=1\}$. In symbols,
$$ 
in_\succ(J_\C) \,\, = \,\, O_\A . 
$$
In fact, in the unimodular case,
Theorem~\ref{th:arrang_resolution} (b)
is precisely  Theorem~4.4 in \cite{BPS}.

\begin{corollary}
The Betti numbers of the unimodular Lawrence ideal $J_\C$,
and all its initial ideals $in_\succ(J_\C)$, are the
coefficients of the cocharacteristic polynomial $\psi_\C$.
\end{corollary}

We close this section with a non-trivial example.
Let $n=9, d=4$ and consider
$$ B^T \quad = \quad
\bordermatrix{ &
 x_{11}&x_{12}&x_{13} & x_{21} & x_{22} & x_{23} & x_{31} & x_{32} & x_{33} \cr
&   1  &  -1  &    0  &    -1  &     1  &     0  &      0 &      0 &      0 \cr
&   0  &   1  &   -1  &      0 &    -1  &     1  &      0 &      0 &      0 \cr
&   0  &   0  &    0  &     1  &    -1  &     0  &     -1 &      1 &      0 \cr
&   0  &   0  &    0  &      0 &     1  &    -1  &      0 &     -1 &    1 \cr}
$$
All nonzero $4 \times 4$-minors of this matrix are $-1$ or $+1$, and 
hence we get a unimodular central arrangement $\C$ of nine
hyperplanes in $\R^4$. This is the {\it cographic arrangement}
associated with the complete bipartite graph $\,K_{3,3} $.
The nine hyperplane variables
$x_{ij}$ represent edges in $ \,K_{3,3} $.
The associated Lawrence ideal can be computed by saturation
(e.g.~in {\tt Macaulay2})
from  (binomials representing) the four rows of $B^T$:
\begin{eqnarray*}
&  J_B  \quad = \quad
 \langle \,\,
x_{11} x_{22}  y_{12} y_{21}  - 
x_{12} x_{21}  y_{11} y_{22} \, , \,
\,x_{12} x_{23}  y_{13} y_{22}  - 
x_{13} x_{22}  y_{12} y_{23}  \, , \, \\ 
& \!\!\! x_{21} x_{32}  y_{22} y_{31}  - 
x_{22} x_{31}  y_{21} y_{32} , 
x_{22} x_{33}  y_{23} y_{32}  - 
x_{23} x_{32}  y_{22} y_{33}  \, \rangle
: ( \prod_{1 \leq i,j \leq 3} x_{ij} y_{ij})^\infty
\end{eqnarray*}
This ideal has $15$ minimal generators, corresponding to the $15$
circuits in the directed graph $K_{3,3}$. 
A typical initial monomial ideal
$\, in_\prec(J_B) = O_\A \,$ looks as follows:
\begin{eqnarray*}
&  \bigl\langle
x_{11} x_{22} y_{12} y_{21} , \,
x_{11} x_{23} y_{13} y_{21} ,\,
x_{11} x_{32} y_{12} y_{31} ,\,
x_{11} x_{33} y_{13} y_{31} ,\,
x_{12} x_{23} y_{13} y_{22} , \\ & 
x_{12} x_{33} y_{13} y_{32} , \,\,
x_{21} x_{32} y_{22} y_{31} , \,\, 
x_{21} x_{33} y_{23} y_{31} , \, \,
x_{22} x_{33} y_{23} y_{32} , \,  \\ &
x_{11} x_{22} x_{33} y_{13} y_{21} y_{32} , \, 
x_{11} x_{22} x_{33} y_{12} y_{23} y_{31} , \,
x_{11} x_{23} x_{32} y_{13} y_{22} y_{31} , \,  \\ &
x_{12} x_{21} x_{33} y_{11} y_{23} y_{32} , \, 
x_{12} x_{21} x_{33} y_{13} y_{22} y_{31} , \, 
x_{13} x_{21} x_{32} y_{12} y_{23} y_{31} \bigr\rangle .
\end{eqnarray*}
This is the oriented matroid ideal of 
the $3$-dimensional affine arrangement $\A$ gotten from $\C$ by 
taking a  vector $w \in \R^4$ with strictly positive coordinates.
This ideal is the intersection of $81$ monomial primes, one
for each spanning tree of $K_{3,3}$. By  Theorem~\ref{lawlaw},
they form a triangulation of a $13$-dimensional Lawrence polytope, 
which is given by its centrally symmetric Gale diagram
$(B^T,-B^T)$ as in \cite[Proposition~9.3.2 (b)]{BLSWZ}.
Resolving this ideal (e.g.~in {\tt Macaulay2}), we obtain
the cocharacteristic polynomial:
\begin{equation}
\label{cochar3by3}
 \psi_\C (q) \quad = \quad
 1   +   15 q +    48 q^2  +  54 q^3 +  20 q^4. 
\end{equation}
It was asked in \cite[\S 5]{BPS} what such Betti numbers
arising from graphic and cographic ideals are in general.
This question will be answered in the following section.

\bigskip

\section{Graphic and Cographic matroids}

Two classes of (central) unimodular arrangements arise from graphs:
{\it graphic\/} and {\it cographic\/} arrangements. Our aim is to
compute their cocharacteristic polynomial. This task is easier for
graphic arrangements which will be treated first.
Cographic arrangements are more challenging and will
be discussed further below.

Fix a connected graph $G$ with vertices $[d] = 
\{1, \ldots, d\}$ and edges $E \subset [d] \times [d]$. Let
$V = \{(v_1, \ldots, v_d)\in \R^d : v_1 + \ldots + v_d =0\}\simeq \R^{d-1}$.
The graphic arrangement $\C_G$ is the  arrangement in $V$
given by the hyperplanes $v_i = v_j$ for  $(i,j)\in E$. It is unimodular
\cite{Wh}.
For each subset $S \subset [d]$ we get an {\it induced subgraph\/}
$\, G|_S \, = \, \bigl( S, \, E \cap (S \times S) \bigr)$.
For a partition $\pi$ of $[d]$, we denote by  $G/\pi$ the
graph obtained from $G$ by contracting all edges whose
vertices lie in the same part of $\pi$. 
The intersection lattice $L_G$ of the graphic arrangement
$\C_G$ has the following
well-known description in terms of the  {\it partition lattice} $\Pi_d$.
See e.g.~\cite{Zas} for proofs and references.

\begin{proposition}
\label{prop:graphic-lattice}
The intersection lattice $L_G$ is isomorphic to the sublattice of the
partition lattice~$\Pi_{d}$ consisting of partitions $\pi$ such 
that, for each part $S$ of $\pi$, the subgraph $G|_S$ is connected.
The element $V_\pi$ of $ L_G$ corresponding to $\pi \in \Pi_d$ 
is the intersection of the hyperplanes $\{v_i = v_j\}$
for pairs $i, j$ in the same part of $\pi$.
The dimension of~$V_\pi$ is equal to the number of parts
of~$\pi$ minus~$1$.  The interval $[V_\pi,\hat 1 \,]$ of the 
intersection lattice $L_G$ is isomorphic to the
intersection lattice $L_{G/\pi}$.
\end{proposition}

We write $\,\mu(G) \, = \, |\mu_{L_G}(\hat{0}, \hat{1})| \,$
for the M\"obius invariant of the intersection lattice.
Thus $\mu(G)$ equals the
{\it Cohen-Macaulay type} (top Betti number) of the 
matroid ideal 
$$ 
M_G \quad = \quad \bigcap \,\bigl\{
\langle  \,x_{ij} \,:\, (i,j) \in F \,\rangle \,\,\, |\,\,\,
F \subseteq E \,\,\,\hbox{is a spanning tree of } \,\, G \,\bigr\}. 
$$
{}From Proposition~\ref{prop:graphic-lattice}  and (\ref{stanley}),
we conclude that all the lower Betti numbers
can be expressed in terms of the M\"obius invariants of the
contractions $G/\pi$ of $G$.

\begin{corollary} The cocharacteristic polynomial of
the graphic arrangement $\C_G$ is
$$
\psi_{\C_G} (q)  \quad = \quad 
\sum_{\pi \in L_G} 
\mu(G/\pi) \cdot q^{|\pi| -1} .
$$
\end{corollary}

This reduces our problem to computing the
M\"obius invariant $\mu(G)$ of a graph $G$.
Green and Zaslavsky~\cite{GZ} found the following 
combinatorial formula.
An {\it orientation\/} of the graph~$G$ is a choice, for each edge 
$(i,j)$ of~$G$, of one of the two possible directions:
$i\to j$ or $j\to i$.  An orientation is {\it acyclic\/}
if there is no directed cycle.

\begin{proposition}
\label{prop:only-source}
Fix a vertex $i$ of $G$. Then $\mu(G)$ equals the number
of acyclic orientations of~$G$ such that,
for any vertex $j$, there is a directed path from $i$ to $j$.
\end{proposition}

\proof
The regions of the graphic arrangement~$\C_G$ are in one-to-one 
correspondence with the acyclic orientations of~$G$:
the region corresponding to an acyclic orientation~$o$ is given by the
inequalities $x_i>x_j$ for any directed edge $i\to j$ in~$o$.

The linear functional  $\,w \, : \, (u_1,\dots,u_d) \mapsto
u_i \, $ is generic for the arrangement $\C_{G}$. The M\"obius invariant
$\mu(G)$ equals the number of regions of $\C_G$ which are bounded
below with respect to $w$. We claim that the acyclic orientations
corresponding to the $w$-bounded regions are precisely the ones
given in our assertion.

Suppose that, for any vertex~$j$ in~$G$, there is a directed path
$i\to  \cdots\to j$.  For any point $(u_1,\dots,u_d)$ 
of the corresponding region, this path implies $u_i>\cdots > u_j$. 
The condition $u_1+\cdots+u_{m}=0$ forces  $w(u) = u_i > 0$.
This implies that the region is $w$-positive.
Conversely, consider any acyclic orientation which does
not satisfy the condition in Proposition \ref{prop:only-source}.
Then there exists a vertex $j\ne i$ 
which is a source of $o$. Then
the vector $v=(-1,\dots,-1,d-1,-1,\dots,-1)$, where
$d-1$ is in the $j$-th coordinate, belongs to the closure of the 
region associated with~$o$.  But $w(v)=-1$.  Hence the region is
not $w$-positive.
\endproof

The above discussion can be translated into a combinatorial
recipe for writing the minimal free resolution of graphic ideals
$M_G$, where each syzygy is indexed by a certain acyclic orientation 
of a graph $G/\pi$. For the case of the {\it complete graph\/} $G=K_{d}$,
we recover the resolution in \cite[Theorem~5.3]{BPS}.
Note that the  intersection lattice $L_{K_{d}}$ is isomorphic to
the partition lattice $\Pi_{d}$.  For any partition $\pi$ 
of $\{1,\ldots,d\}$ with
$i+1$ parts, $K_d/\pi$ is isomorphic to $K_{i+1}$. The number of
such partitions equals  $S(d,i+1)$, the {\it Stirling number} of 
the second kind. The number of acyclic orientations of $K_{i+1}$ 
with a unique fixed source  equals $i\,! $. We conclude

\begin{corollary}
The number of minimal $i$-th syzygies of
$M_{K_d}$ equals $i!\ S(d,i+1)$.
\end{corollary}

\begin{remark} \label{foata} \rm Vic Reiner suggested to us the following
combinatorial interpretation of $\mu(G)$.
It can be derived from Proposition \ref{prop:only-source}.
For any graph $G$, the M\"obius invariant
$\mu(G)$ counts the number of equivalence classes of linear
orderings of the vertices of $G$, under the equivalence
relation generated by the following operations:
\begin{enumerate}
\item[$\bullet$] commuting two adjacent vertices $v,v'$
in the ordering if $\{v,v'\}$ is not an edge of $G$,
\item[$\bullet$] cyclically shifting the entire order, i.e. 
$\,\,v_1 v_2 \ldots v_n \leftrightarrow v_2 \ldots v_n v_1 $.
\end{enumerate}
Invariance under the second operation makes this interpretation
convenient for writing down the minimal free resolution of the 
graphic Lawrence  ideals in \cite[\S 5]{BPS}.

Another application arises when $(W,S)$ is a Coxeter system 
and $G$ its Coxeter graph (considered without its edge labels).
Suppose $S=\{s_1,\ldots,s_n\}$. Then $\mu(G)$ 
counts the number of Coxeter elements $s_{i_1} \ldots s_{i_n}$
of $G$ up to the equivalence relation $\,
s_{i_1} s_{i_2} \ldots s_{i_n} \leftrightarrow 
s_{i_2} \ldots s_{i_n} s_{i_1} $.
\end{remark}
\noindent

\medskip

We now come to the cographic arrangement $\C^\perp_G$,
whose matroid is dual to that of  $\C_G$. 
Fix a directed graph $G$ on $[d]$ with edges $E$,
where $G$ is allowed to have loops and multiple edges.
We associate with $G$ the multiset
of vectors $\{v_e \in \Z^{d} : e\in E\}$, where 
for an edge $e=(i\to j)$,
the $i$th coordinate of $v_e$ is $1$, the $j$th coordinate
is $-1$, and all other coordinates are~$0$.
Set $v_e=0$ for a loop $e=(i\to i)$ of $G$.
Let $V_G=\{\lambda: E \rightarrow \R \mid \sum_{e\in E} \lambda(e)v_e=0\}$.
Note that $V_G$ is a vector space of dimension
$\, \#\{\textrm{edges}\}-\#\{\textrm{vertices}\} + 
\#\{\textrm{connected components}\}$.
The {\it cographic arrangement\/} $\C^\perp_G$ is the arrangement 
in $V_G$ given by hyperplanes
$ H_e=\{\lambda\in V_G : \lambda(e)=0\}$ for $e\in E$.
It is unimodular \cite{Wh}.
We write $\,\mu^\perp(G) \, = \, |\mu_{L^\perp_G}(\hat{0}, \hat{1})| \,$
for the M\"obius invariant of the intersection lattice $L^\perp_G$
of $\C^\perp_G$, and we refer to this number as the
{\it M\"obius coinvariant\/} of~$G$.
Thus $\mu^\perp(G)$ is the Cohen-Macaulay type of the
cographic ideal $J_{C^\perp_G}$ in \cite[\S 5]{BPS}.

\smallskip

\begin{remark} {\rm
The characteristic polynomial of a matroid can 
be expressed via the Tutte
dichromatic polynomial~\cite{T}.  Thus M\"obius invariant and coinvariant of a
graph $G$ are certain values of the Tutte polynomial: $\mu(G)=T_G(1,0)$ and
$\mu^\perp(G)=T_G(0,1)$.  We do not know, however, how to express the
cocharacteristic polynomial $\psi(q)$ in terms of the Tutte polynomial.
}
\end{remark}

A formula for the Tutte polynomial due to 
Gessel and Sagan~\cite[Theorem~2.1]{GS} implies:

\begin{proposition} 
\label{Coinvariant}
The M\"obius coinvariant of $G$ is
$\,
\mu^\perp(G)  \,=\, \sum_{F\subseteq G} (-1)^{d-|F|-1},
$
where the sum is over all forests in $G$
and $|F|$ denotes the number of edges in $F$.
\end{proposition}

We shall derive explicit formulas for the 
M\"obius coinvariant of complete  and complete bipartite graphs.
A subgraph~$M$ of a graph $G$ is called a {\it partial matching\/}
if it is a collection of pairwise disjoint edges of the graph.
For a partial matching~$M$, let $a(M)$ be the number of vertices 
of $G$ that have degree~$0$ in~$M$.  The {\it Hermite polynomial\/}
$H_n(x)$, $n\geq 0$, is the generating function of partial matchings
in the complete graph $K_n$:
$$
H_n(x) = \sum_M x^{a(M)},
$$
where the sum is over all partial matchings in~$K_n$.  In particular,
$H_0(x)=1$.  Set also $H_{-1}(x)=0$.
The main result of this section is the following formula:

\begin{theorem}  \label{th:mu(K_m)}
The M\"obius coinvariant of the complete graph $K_m$ equals
\begin{equation}  \label{eq:mu_m}
\mu^\perp(K_m)  \quad = \quad 
(m-2)\, H_{m-3}(m-1)\,,\qquad m\geq 2\,.
\end{equation}
\end{theorem}

A few initial numbers $\mu^\perp(K_m)$ are given below.
$$
\begin{array}{lccccccccccc}
m      & \quad &  2 & 3 & 4 & 5  & 6   & 7    & 8      & 9    & 
10    & \cdots 
\\[.05in]
\mu^\perp(K_m)  &       &  0 & 1 & 6 & 51 & 560 & 7575 & 122052 & 2285353 &
48803904  & \cdots
\end{array}
$$
The proof of Theorem~\ref{th:mu(K_m)} relies on several auxiliary results
and will be given below. The next proposition summarizes well-known  
properties of Hermite polynomials.

\begin{proposition}
\label{prop:Hermite}
The Hermite polynomial $H_n(x)$ satisfies  the recurrence
\begin{equation}
\label{eq:H-recurrence}
\begin{array}{l}
H_{-1}(x)=0,\quad H_0(x)=1,\\[.05in]
H_{n+1}(x) = x\, H_n(x) + n\, H_{n-1}(x),\quad n\geq 0.
\end{array}
\end{equation}
It is given explicitly by the formula
 $$ H_n(x) \quad =\quad x^n + \sum_{k\geq1}^{[n/2]} \binom{n}{2k} (2k-1)!!\, 
x^{n-2k}, $$
where $(2k-1)!!=(2k-1)(2k-3)(2k-5)\cdots 3\cdot 1$.
\end{proposition}

\proof
In a partial matching the first vertex has either degree~$0$
or~$1$.  This gives two terms in the right-hand side 
of the recurrence~(\ref{eq:H-recurrence}).  The formula for $H_n(x)$
follows from the fact that there are $(2k-1)!!$ matchings with $k$ edges 
on $2k$ vertices.
\endproof

Returning to general cographic arrangements,
recall that an edge $e$ of the graph~$G$ is called an {\it isthmus\/}
if $G \backslash e$ has more connected components than~$G$;
a graph is called {\it isthmus-free\/} if no edge of~$G$ is an isthmus.
The minimal nonempty isthmus-free subgraphs of~$G$
are the {\it cycles\/} of~$G$.
For a subgraph~$H$ of~$G$, denote by $G/H$  the graph obtained 
from~$G$ by contracting the edges of~$H$.  Note that  $G/H$ may
have loops and multiple edges even if $G$ does not.
The following result appears in \cite{GZ}.

\begin{proposition} 
\label{prop:cographic}
The intersection lattice  $L^\perp_G$ of the cographic arrangement
is isomorphic to the lattice of isthmus-free subgraphs of~$G$
ordered by reverse inclusion.
The element of the intersection lattice that corresponds to
an isthmus-free subgraph~$H$ is $V_H\subset V_G$.
The coatoms of the lattice $L^\perp_G$ are  the cycles of~$G$.
 For two isthmus-free subgraphs~$H\supset K$ of~$G$, the interval
$[V_H,V_K]$ of the intersection lattice  $L^\perp_G$  
is isomorphic to the interval $[\hat 0,\hat 1]$ 
of the intersection lattice $L^\perp_{H/K}$.
\end{proposition}

Proposition~\ref{prop:psi-recurrence} implies the following recurrence 
for the cocharacteristic polynomial $\psi_{\C^\perp_G}(q)$
of the cographic arrangement $\C^\perp_G$.

\begin{corollary}
Let $e$ be an edge of the graph~$G$. Then
\begin{equation}
\label{eq:psi_recurr}
\psi_{\C^\perp_G}(q) = \psi_{{\C^\perp_{G\setminus e}}}(q) + 
                     q \sum_C \psi_{\C^\perp_{G/C}}(q),
\end{equation}
where the sum is over all cycles~$C$ of $G$ that contain $e$.
\end{corollary}

Considering  terms of the highest degree in (\ref{eq:psi_recurr}),
we obtain 

\begin{corollary}
\label{cor:mu_recurr}
If $e$ is any edge of $G$ that is not an isthmus, then
\begin{equation}
\label{eq:mu_recurr}
\mu^\perp(G) \quad = \quad  \sum_C \mu^\perp(G/C),
\end{equation}
where the sum is over all cycles~$C$ of $G$ that contain $e$.
\end{corollary}

Note that $\mu^\perp(G)$ is equal to the M\"obius coinvariant of
the graph ${\widetilde G}$ obtained from~$G$
by removing all loops and isthmuses.
Thus when we use  relation~(\ref{eq:mu_recurr}) to calculate
$\mu^\perp(G)$, we may remove all new loops obtained 
after contracting the cycle~$C$.

 We are ready to prove Theorem~\ref{th:mu(K_m)}.
For $n\geq 0$ and $k\geq 1$, define  $K_n^{(k)}$  to be 
 the complete graph $K_n$ on the vertices $1,\dots,n$,
together with one additional vertex $n+1$ (root) connected  
to each vertex $1,\dots,n$ by $k$ edges.
Let $\mu_n^{(k)}=\mu^\perp(K_n^{(k)})$ be the M\"obius-coinvariant
of the graph~$K_n^{(k)}$.
Note that $K_m = K_{m-1}^{(1)}$ and $\mu^\perp(K_m)=\mu_{m-1}^{(1)}$.
Theorem~\ref{th:mu(K_m)} can be extended as follows:

\begin{proposition}
\label{prop:mu_nk}
$
\mu_{n}^{(k)}= H_n(n+k-1) - n\, H_{n-1}(n+k-1)\,
$
for $n,k\geq 1$.
\end{proposition}

\proof 
We  utilize Corollary~\ref{cor:mu_recurr}. 
Select an edge $e=(n,n+1)$ of the graph $K_n^{(k)}$.
There are $k-1$ choices for a cycle~$C$ of length~$2$ that contains
the edge~$e$, and the graph $K_n^{(k)}/C$, after removing 
 loops, is isomorphic to $K_{n-1}^{(k+1)}$.
There are $(n-1)\,k$ choices for a cycle~$C$ of length~$3$ that contains
the edge~$e$, and the graph $K_n^{(k)}/C$, after removing loops,
is isomorphic to $K_{n-2}^{(k+2)}$.
%For cycles of length $4$, we have $(n-1)(n-2)\,k$ choices, and
%obtain a graph that is isomorphic to $K_{n-3}^{(k+3)}$, etc.
In general, for cycles of length $l\geq 3$, there are
$k\,(n-1)(n-2)\cdots(n-l+2)$ choices, and we obtain a graph that is
isomorphic to $K_{n-l+1}^{(k+l-1)}$.
Equation~(\ref{eq:mu_recurr}) implies the following
recurrence for  $\mu_n^{(k)}$:
\begin{equation}
\label{eq:mu_nk-recur}
\begin{array}{l}
\mu_n^{(k)} \quad = \quad (k-1)\, \mu_{n-1}^{(k+1)} 
+ k\, (n-1)\,  \mu_{n-2}^{(k+2)} + \\[.05in]
+ k\, (n-1)(n-2)\,  \mu_{n-3}^{(k+3)} + k\, (n-1)(n-2)(n-3)\,
\mu_{n-4}^{(k+4)} + \cdots \,,
\end{array}
\end{equation}
which, together with the initial condition 
$\mu_0^{(k)}=1$, defines the numbers $\mu_n^{(k)}$ uniquely.
Set
$$
b_n^{(k)} \quad = \quad \mu_n^{(k)} + n \, \mu_{n-1}^{(k+1)} +
n(n-1) \, \mu_{n-2}^{(k+2)} +\cdots+
n(n-1)\cdots 1 \, \mu_{0}^{(k+n)}.
$$
Then $\mu_n^{(k)} = b_n^{(k)}  -n b_{n-1}^{(k+1)}$ and 
the relation~(\ref{eq:mu_nk-recur}) can be rewritten as
$$
b_n^{(k)} -n\, b_{n-1}^{(k+1)} = (k-1)
\left(b_{n-1}^{(k+1)} -(n-1)\, b_{n-2}^{(k+2)}\right) +
k\,(n-1)\,b_{n-2}^{(k+2)},
$$
or, simplifying, as
\begin{equation}
\label{eq:bnk}
b_n^{(k)} = (n+k-1) b_{n-1}^{(k+1)} +(n-1)\,b_{n-2}^{(k+2)}.
\end{equation}
We claim that $b_n^{(k)} = H_n(n+k-1)$.  Indeed, $b_0^{(k)}=1$, 
$b_1{(k)}=k$, and equation~(\ref{eq:bnk}) is equivalent to 
the defining relation~(\ref{eq:H-recurrence}) for the Hermite polynomials.
Hence $\mu_n^{(k)} = b_n^{(k)}  -n b_{n-1}^{(k+1)} = 
H_n(n+k-1) - n H_{n-1}(n+k-1)$.
\endproof

\proofof{Theorem~\ref{th:mu(K_m)}}
By Proposition~\ref{prop:mu_nk} and equation~(\ref{eq:H-recurrence}), 
$$
\mu^\perp(K_m) = \mu_{m-1}^{(1)} = H_{m-1}(m-1) - (m-1)\,H_{m-2}(m-1) =
(m-2)\,H_{m-3}(m-1).
$$
\endproof

We now discuss a bipartite analog of Hermite polynomials.
For a partial matching $M$ in the complete bipartite graph $K_{m,n}$,
denote by $a(M)$ the number of  vertices in the first 
part that have degree 0 in $M$, and by $b(M)$ the number of vertices
in the second part that have degree~$0$.
Define
$$
H_{m,n}(x,y) \quad = \quad \sum_{M} x^{a(M)}\, y^{b(M)}\,,
$$
where the sum is over all partial matchings in $K_{m,n}$.
In particular $H_{m,0}=x^m$ and $H_{0,n}=y^n$.
Set also $H_{m,-1} = H_{-1,n}=0$.
The following statement is a bipartite analogue of Theorem~\ref{th:mu(K_m)}.

\begin{theorem}
\label{th:mu-bipartite}
The M\"obius coinvariant of the complete bipartite graph $K_{m,n}$
 equals 
$$
\mu^\perp(K_{m,n}) \quad = \quad (m-1)(n-1) H_{m-2,\,n-2}(n-1,m-1)\,,
\qquad m,n\geq 1\,.
$$
\end{theorem}

\smallskip

Analogously to Proposition~\ref{prop:Hermite}
we have

\begin{proposition} 
The polynomial $H_{m,n}(x,y)$ is given by
$$
H_{m,n}(x,y) = \sum_{k=0}^{\min(m,n)}
\binom{m}{ k} \binom{n}{ k}k!\, x^{m-k} y^{n-k}\,.
$$
It  satisfies the following recurrence relations:
\begin{equation}
\label{eq:relations-bipartite}
\begin{array}{l}
H_{m,n}(x,y) = x\,H_{m-1,n}(x,y) + n\,H_{m-1,n-1} (x,y), \\[.1in]
H_{m,n}(x,y) = y\,H_{m,n-1}(x,y) + m\,H_{m-1,n-1} (x,y), \\[.1in]
H_{m,0}=x^m,\qquad H_{0,n}=y^n.
\end{array}
\end{equation}
\end{proposition}

\proof 
The first formula is obtained by counting the partial matchings in $K_{m,n}$.
The recurrence relations~(\ref{eq:relations-bipartite}) are obtained
by distinguishing two cases when the first vertex in the first (second) part 
of $K_{m,n}$ has degree $0$ or $1$ in a partial matching.
\endproof

Let us define the graph $K_{m,n}^{(k,l)}$ as the complete 
bipartite graph $K_{m,n}$ with an additional vertex~$v$
such that $v$ is connected by $k$ edges with each vertex in
the first part and by $l$ edges with each vertex in the second part.
Let $\mu_{m,n}^{(k,l)}=\mu^{\perp}(K_{m,n}^{(k,l)})$ be the M\"obius
coinvariant of this graph.  Note that $K_{m,n}=K_{m,n-1}^{(1,0)}$ and,
thus, $\mu^\perp(K_{m,n})=\mu_{m,n-1}^{(1,0)}$.
Theorem~\ref{th:mu-bipartite} can be extended as follows:

\begin{proposition}
\label{prop:mu_mnkl}
We have
$$
\mu_{m,n}^{(k,l)}=H_{m,n}(n+k-1,m+l-1) - 
m\,n\ H_{m-1,\,n-1}(n+k-1,m+l-1).
$$
\end{proposition}

\proof
Our proof is similar to that of Proposition~\ref{prop:mu_nk}.
We  utilize Corollary~\ref{cor:mu_recurr}. 
Select an edge $e$ of the graph $K_{m,n}^{(k,l)}$ that joins the additional
vertex $v$ with a vertex from the first part.
There are $k-1$ choices for a cycle~$C$ of length~$2$ that contains
the edge~$e$, and the graph $K_{m,n}^{(k,l)}/C$, after removing 
 loops, is isomorphic to $K_{m-1,n}^{(k,l+1)}$.
There are $n\,l$ choices for a cycle~$C$ of length~$3$ that contains
the edge~$e$, and the graph $K_{m,n}^{(k,l)}/C$, after removing loops,
is isomorphic to $K_{m-1,n-1}^{(k+1,l+1)}$.
For cycles of length $4$, we have $n(m-1)\,k$ choices, and
obtain a graph isomorphic to $K_{m-2,n-1}^{(k+1,l+2)}$, etc.
In general,  for cycles of odd length $2r+1\geq 3$, we have 
$l\,n(m-1)(n-1)(m-2) \cdots (m-r+1)(n-r+1)$ choices, and we obtain a graph
isomorphic to $K_{m-r,\,n-r}^{(k+r,\,l+r)}$.
For cycles of even length $2r+2\geq 4$, we have 
$k\,n(m-1)(n-1)(m-2)\cdots (n-r+1)(m-r)$ choices, and we obtain a graph 
isomorphic to $K_{m-r-1,\,n-r}^{(k+r,\,l+r+1)}$.
The equation~(\ref{eq:mu_recurr}) implies the following
recurrence for  $\mu_{m,n}^{(k,l)}$:
\begin{equation}
\label{eq:mu_mnkl-recur}
\begin{array}{l}
\mu_{m,n}^{(k,l)} = (k-1)\, \mu_{m-1,n}^{(k,l+1)} 
+ l\, n\,  \mu_{m-1,n-1}^{(k+1,l+1)} + 
k\, n(m-1)\,  \mu_{m-2,n-1}^{(k+1,l+2)} + \\[.05in] 
l\, n(m-1)(n-1)\, \mu_{m-2,n-2}^{(k+2,l+2)} 
+ k\, n(m-1)(n-1)(m-2)\,  \mu_{m-3,n-2}^{(k+2,l+3)} + \cdots \,,
\end{array}
\end{equation}
which, together with the initial conditions 
$\mu_{0,n}^{(k,l)}=(l-1)^n$ and 
$\mu_{m,0}^{(k,l)}=(k-1)^m$, unambiguously defines the numbers 
$\mu_{m,n}^{(k,l)}$.
Let us fix the numbers $p=k+n-1$ and $q=l+m-1$ and write 
$\mu_{m,n}$ for $\mu_{m,n}^{(p-n+1,q-m+1)}$.
Set
$$
b_{m,n} =\mu_{m,n}+ n\, m\, \mu_{m-1,n-1} +
n(n-1)\,m(m-1) \, \mu_{m-2,n-2}+\cdots.
$$
Then $\mu_{m,n} = b_{m,n} -m\,n b_{m-1,n-1}$ and 
the relation~(\ref{eq:mu_mnkl-recur}) can be rewritten as
$$
\begin{array}{l}
b_{m,n} -m\,n\, b_{m-1,n-1} 
=-\left(b_{m-1,n} -(m-1)\,n\, b_{m-2,n-1}\right)+ \\[.05in]
\qquad\qquad\qquad\qquad
+ (p-n+1) \, b_{m-1,n} + (q-m+1)\,n\, b_{m-1,n-1}
\end{array}
$$
or, simplifying, as
\begin{equation}
\label{eq:bnm-bipartite}
b_{m,n} = (p-n)\, b_{m-1,n} +(q+1)n\,b_{m-1,n-1} + (m-1)n\, b_{m-2,n-1}.
\end{equation}
This relation, together with the initial conditions $b_{0,n}=q^n$,
$b_{m,0}=p^m$, $b_{-1,n}=b_{m,-1}=0$, uniquely determines the numbers
$b_{m,n}$.

We claim that $b_{m,n} = H_{m,n}(p,q)$.  Indeed, the above initial conditions
are satisfied by $H_{m,n}(p,q)$ and~(\ref{eq:bnm-bipartite}) follows from
the defining relations~(\ref{eq:relations-bipartite}) for the bipartite
Hermite polynomials.  In order to see this,
we write by~(\ref{eq:relations-bipartite})
$$
\begin{array}{c}
H_{m,n}(p,q)=p\,H_{m-1,\,n}(p,q)+n\,H_{m-1,\,n-1}(p,q),\\[.05in]
n\,H_{m-1,\,n}(p,q)=n\,q\,H_{m-1,\,n-1}(p,q)+n(m-1)\,H_{m-2,\,n-1}(p,q)\,.
\end{array}
$$
The sum of these two equations is equivalent to the 
equation~(\ref{eq:bnm-bipartite}).
Hence $\mu_{m,n}^{(k,l)} = b_{m,n}  - m\,n\, b_{m-1,n-1}=
H_{m,n}(p,q) - m\,n\, H_{m-1,\,n-1}(p,q)$.
\endproof

An alternative expression for $\mu_{m,n}^{(k,l)}$ can be deduced
from Proposition~\ref{prop:mu_mnkl}:
\begin{equation}
\label{eq:alt-mu_mnkl}
\mu_{m,n}^{(k,l)} = \sum_{r=0}^{\min(m,n)}
(1-r)\,\binom{m}{ r} \binom{n}{ r} r!\, (n+k-1)^{m-r} (m+l-1)^{n-r}\,.
\end{equation}

\proofof{Theorem~\ref{th:mu-bipartite}}
By Proposition~\ref{prop:mu_mnkl} and the recurrence 
relations~(\ref{eq:relations-bipartite}),
$$
 \begin{array}{l}
\mu^\perp(K_{m,n}) \\[.05in] 
\ = \,\, \mu_{m,n-1}^{(1,0)} 
 \,\, =\,\, H_{m,n-1}(n-1,m-1) - 
m(n-1)\,H_{m-1,\,n-2}(n-1,m-1)  \\[.05in]
\  = \,\, (n-1)\,H_{m-1,\,n-1}(n-1,m-1)-(m-1)(n-1)\,H_{m-1,\,n-2}(n-1,m-1) 
\\[.05in]
\ = \,\,\, (m-1)(n-1)\,H_{m-2,\,n-2}(n-1,m-1)\,.
\end{array}
$$
\endproof

\smallskip

For a commutative algebra example illustrating
Theorem~\ref{th:mu-bipartite}, consider the Lawrence ideal
$\, J_B \subset \k[
x_{11}, \ldots, x_{33}, 
y_{11}, \ldots, y_{33}]\,$
associated with the bipartite graph $K_{3,3}$.
This is the Lawrence lifting of  the ideal of $2 \times 2$-minors
of a  generic $3 \times 3$-matrix. It
was discussed in the end of Section 4.
Its Cohen-Macaulay type is
$$ \mu^\perp(K_{3,3}) \quad = \quad
(3-1)\cdot (3-1) \cdot H_{1,1}(2,2) 
\quad = \quad
2 \cdot 2 \cdot 5 \quad = \quad 20 . $$
This is the leading coefficient of the cocharacteristic
polynomial in equation (\ref{cochar3by3}).

\bigskip \bigskip

\noindent {\bf Acknowledgements}:
We wish to thank Vic Reiner and G\"unter Ziegler
for valuable communications. Their ideas and suggestions
have been incoorporated
in Proposition \ref{contractible},
Theorem~\ref{f-vector} and Remark \ref{foata}.
Alexander Postnikov was partially supported by NSF grant
\#DMS-9840383. Bernd Sturmfels was partially supported by NSF grant
    \#DMS-9970254.

\bigskip \bigskip

\end{document}